\documentclass[MSNbibl,number,citesort,seceqn,dvips]{arxbj}
\usepackage{esint}

\aid{0}
\volume{19}
\issue{4}
\pubyear{2013}
\firstpage{1419}
\lastpage{1448}
\doi{10.3150/12-BEJSP08} 

\makeatletter
\newcommand{\RV}{\mathrm{RV}}
\newcommand{\cpt}{\mathcal{K}}
\newcommand{\E}{\mathbb{E}}
\newcommand{\cbfn}{\mathcal{C}}
\newcommand{\mplus}{\mathbb{M}_{+}}
\newcommand{\Mp}{\mathbb{M}_{p}}
\newcommand{\leb}{\mathbb{LEB}}
\newcommand{\prm}{\operatorname{PRM}}
\newcommand{\pp}[1]{\epsilon_{#1}}
\newcommand{\indfn}{\mathbf{1}}
\newcommand{\wc}{\Rightarrow}
\newcommand{\bX}{\boldsymbol{X}}
\newcommand{\given}{\mid}
\newcommand{\EP}{\mathsf{E}}
\renewcommand{\P}{\mathsf{P}}
\newcommand{\eqref}[1]{(\ref{#1})}
\newcommand{\bxi}{\boldsymbol{\xi}}
\newcommand{\bT}{\boldsymbol{T}}
\newcommand{\conv}{\longrightarrow}
\newcommand{\ind}[1]{\mathbf{1}_{\{#1\}}}
\newcommand{\smax}{\vee}
\newcommand{\R}{\mathbb{R}}
\newcommand{\rvect}[1]{\boldsymbol{#1}}
\newcommand{\bz}{\boldsymbol{0}}
\newcommand{\bdry}{\partial}
\newcommand{\smin}{\wedge}

\newcommand{\eb}{y}
\newcommand{\tat}{\tau(t)}

\newcommand{\tbn}{\tau(b_n)}
\newcommand{\limpp}{\zeta}
\newcommand{\nua}{\nu_\alpha}
\newcommand{\tk}{\tilde{\tau}_{k+1}(b_n)}
\newcommand{\tkit}{\tilde{\tau}_k(t)}
\newcommand{\tkoit}{\tilde{\tau}_{k+1}(t)}
\newcommand{\sk}{\sigma_{k+1}(b_n)}
\newcommand{\skt}{\sigma_{k+1}(t)}
\newcommand{\Xkj}{X_{kj}}
\newcommand{\Xko}{X_{k0}}
\newcommand{\Xkm}{X_{km}}
\newcommand{\xikj}{\xi_{k}(j)}
\newcommand{\xiko}{\xi_{k}(0)}
\newcommand{\xikm}{\xi_{k}(m)}
\newcommand{\xivkm}{\bxi_{k,m}}
\newcommand{\ppnveco}{\vartheta_n}
\newcommand{\ppveco}{\vartheta}
\newcommand{\Ld}{\Lambda_\delta}
\newcommand{\ppnvst}{\hat{\eta}_n}
\newcommand{\ppvst}{\hat{\eta}}
\newcommand{\ppnst}{\tilde{\eta}_n}
\newcommand{\ppst}{\tilde{\eta}}
\newcommand{\ppnmain}{\eta_n}
\newcommand{\ppmain}{\eta}
\newcommand{\ta}{\tau_A}
\newcommand{\tkt}{\tau_k(t)}
\newcommand{\tkob}{\tau_{k+1}(b_n)}
\newcommand{\tkb}{\tau_k(b_n)}
\newcommand{\ppmcd}{\chi_n}
\newcommand{\Nn}{N_n}
\newcommand{\Nnd}{\tilde{N}_n}
\newcommand{\bzero}{\boldsymbol0}
\newcommand{\bV}{\rvect{V}}

\newtheorem{them}{Theorem}[section]
\newtheorem{lem}{Lemma}[section]
\newtheorem{prop}{Proposition}[section]
\newremark{rmk}{Remark}
\makeatother

\begin{document}
\begin{frontmatter}

\title{Clustering of Markov chain exceedances}
\runtitle{Markov chain exceedances}

\begin{aug}
\author[1]{\fnms{Sidney I.} \snm{Resnick}\corref{}\thanksref{1}\ead[label=e1]{sir1@cornell.edu}} \and
\author[2]{\fnms{David} \snm{Zeber}\thanksref{2}\ead[label=e2]{dsz5@cornell.edu}}
\runauthor{S.I. Resnick and D. Zeber} 
\address[1]{School of ORIE, Rhodes Hall 284, Cornell University,
Ithaca, NY 14853, USA.\\ \printead{e1}}
\address[2]{Department of Statistical Science, 301 Malott Hall,
Cornell University, Ithaca, NY 14853, USA.\\ \printead{e2}}
\end{aug}


%
\begin{abstract}
The tail chain of a Markov chain can be used to model the dependence
between extreme observations.
For a positive recurrent Markov chain, the tail chain aids in
describing the limit of a sequence of point processes $\{N_n, n \geq1\}$,
consisting of normalized observations plotted against scaled
time points.
Under fairly general conditions on extremal behaviour, $\{N_n\}$
converges to a cluster Poisson process. Our technique
decomposes the sample path of the chain into i.i.d. regenerative cycles
rather than using blocking argument typically employed in the context
of stationarity with mixing.
\end{abstract}


\end{frontmatter}

\section{Introduction}

One of the effects of dependence in a time series is that extremes
tend to cluster. This has applied implications to risk
contagion over time but is also mathematically interesting and the
challenge is to precisely relate the dependence structure to the
clustering.
For Markov dependence,
how do we
describe exceedance clusters?

Point processes powerfully describe extremal behaviour
of certain time series.
Under appropriate conditions on marginal distributions and rapid decay
of dependence as a
function of time lag
for the process $\{X_j :
j\geq0\}$, the exceedance point process $\Nn$ defined by
%
%
\begin{equation}
\label{eqppnn} \Nn \bigl([0,s] \times(a,\infty] \bigr) = \# \{ j \leq sn : X_j
> ab_n \}
\end{equation}
converges weakly to a Poisson limit as $n\rightarrow\infty$, where $b_n
\rightarrow\infty$ is a threshold sequence.
This leads to a number of results on asymptotic distributions of large
order statistics and exceedances of an extreme level. Such results have
been developed in a variety of contexts by
\cite
{leadbetter1983extremes,hsing1988exceedance,balan2009convergence,hsing1987characterization,hsing1989extreme,borkovec2000extremal,novak2002multilevel}.
More specific results exist for regularly varying processes
\cite{basrak2009regularly,davis1995point}, regenerative processes
\cite{asmussen1998extreme,rootzen1988maxima}, and Markov chains
\cite{perfekt1994extremal,yun1998extremal}. Distributions of
functionals of such point processes have been considered in
\cite{yun2000distributions,segers2003functionals,segers2005approximate}.

For stationary processes, the dependence structure causes extremes to
occur in clusters.
The clustering is often summarized using the extremal index
$\theta$ introduced by Leadbetter et
al. \cite{leadbetter1983extremes}, which is related to the asymptotic
mean cluster size.
To obtain a point process convergence result, authors often employ
the big block/little block technique and
mixing conditions, such as Leadbetter's $D(u_n)$
(see \cite{leadbetter1983extremes}), to split the process into approximately
independent and
identically distributed blocks.
With an appropriate choice of block size, extremes within one such
block belong asymptotically to the same cluster.
Under an assumption controlling the extremal behaviour within each
block, such as via the distribution of the number of exceedances, $\Nn
$ generally converges to a limiting compound Poisson process, where the
compounding at each time point approximates the clustering within each block.
For Markov chains, the \emph{tail chain} is an asymptotic process that
models behaviour upon reaching an extreme state; see
\cite
{smith1992extremal,perfekt1994extremal,perfekt1997extreme,segers2007multivariate,resnick2011asymptotics,smith1997markov}.
Point process results for stationary Markov chains employ the tail
chain to specify the compounding in the limit process.
Under Markov dependence, the within-block behaviour is determined
merely by conditions on the marginal distribution and the transition
kernel.
Basrak and Segers \cite{basrak2009regularly} extended
the tail chain model to general multivariate regularly varying stationary
processes.

Rootz\'en \cite{rootzen1988maxima} focuses on regenerative processes,
which split naturally into \textit{cycles}.
In this case, the within-block condition is replaced by an assumption
on the extremal behaviour over a cycle.
The main difference is that the cycles are of random but finite length,
whereas the block size increases deterministically with $n$.
In particular, Rootz\'en shows that convergence of the sequence of
processes counting the number of exceedances
depends on
the asymptotics of
the distribution of the cycle maximum
as well as the marginal distribution.

We combine these two approaches to derive the weak limit of $\{\Nn\}$
when $\{X_n\}$ is a
positive recurrent
Markov chain.
Such chains display a regenerative structure and
in the limit, $N_n$ is approximated by a process consisting of
clusters of points stacked above common time points, each
corresponding to a separate regenerative cycle.
The heights of the points in each cluster are determined by an
independent run of the tail chain.
This paper requires some distributional results for the tail chain
process that were
derived in \cite{resnick2011asymptotics} and Section~\ref{secextrmch}
offers a summary of necessary facts.
We focus on the case of
heavy-tailed marginals, but believe our approach could be extended to
accommodate more general marginal distributions.

\subsection{Notation and conventions}

We review notation and relevant concepts.
In general, bold symbols represent vectors or sequences
and for $\boldsymbol x=(x_1,x_2,\ldots)$, write
$\boldsymbol x_m := (x_1,\ldots,x_m)$.

\begin{tabular}{ll}
$f^{\leftarrow}$ & the left-continuous inverse of a monotone function $f$,
i.e., \\
& $f^{\leftarrow}(x) = \inf\{y : f(y) \geq x\}$. \\
$\RV_{\rho}$ & the class of regularly varying functions with index
$\rho$. \\
$D[0,\infty)$ & the space of real-valued c\`adl\`ag functions on
$[0,\infty)$ endowed with the \\
&Skorohod topology.\\
$D_{\mathrm{left}} [0,\infty)$ & left continuous functions on
$[0,\infty)$ with finite right hand limits and \\
&metrized by the
Skorohod metric.
\end{tabular}

\begin{tabular}{ll}
$D^{\uparrow}[0,\infty)$ & the subspace of $D_{\mathrm
{left}}[0,\infty)$
consisting of non-decreasing functions $f$ with \\
& $f(0)=0$ and $\lim_{x\rightarrow\infty} f(x)=\infty$. \\
$\cpt(\E)$ & the collection of compact subsets of $\E$. \\
$\cbfn(\E)$ & the space of real-valued continuous, bounded functions
on $\E$. \\
$\mathcal{C}^{+}_K(\E)$ & the space of non-negative continuous
functions on $\E$
with compact support. \\
$\mplus(\E)$ & the space of non-negative Radon measures on $\E$.\\
$\Mp(\E)$ & the space of Radon point measures on $\E$. \\
$\leb$ & Lebesgue measure on $\mathbb{R}$. \\
$\prm(\mu)$ & Poisson random measure on $\E$ with mean measure $\mu
$. \\
$\pp{x}(\cdot)$ & point measure at $x$, i.e., $\pp{x}(A) = \indfn_A(x)$. \\
$\nua$ & a measure on $(0,\infty]$ given by $\nua(x,\infty] =
{x}^{-\alpha}$ for $x>0$, $\alpha>0$. \\
$\wc$ & weak convergence of probability measures \cite{billingsley1999convergence}.
\end{tabular}

For a space $\E$ which is locally compact with countable base (for
example, a subset of $[-\infty,\infty]^d$),
a sequence of measures $\{ \mu_n\} \subset\mplus(\E)$ converges
vaguely to $\mu\in\mplus(\E)$ (written $\mu_n\,\mathop
{\longrightarrow}\limits^v\,\mu$) if
$\int_\E f \,d\mu_n \rightarrow\int_\E f\, d\mu$ as $n\rightarrow
\infty$ for
any $f\in\mathcal{C}^{+}_K(\E)$.
The vague topology on $\mplus(\E)$ is metrizable by the vague metric,
$d_v$, i.e., $d_v(\mu_n,\mu) \rightarrow0$ iff $\mu_n \,\mathop
{\rightarrow}\limits^v\,\mu
$. See
\cite{resnick1987extreme,kallenberg1976random,neveu1977processus} for
further details.
A distribution $F$ on $[0,\infty)$ has a regularly varying tail with
index $\alpha>0$, denoted $1-F \in
\RV_{-\alpha}$, if there exists $b(t)\rightarrow\infty$
such that
\[
tF \bigl(b(t)\cdot \bigr) \,\mathop{\longrightarrow}^v\,
\nu_\alpha(\cdot) \qquad\mbox{in } \mplus(0,\infty] \qquad\mbox {as } t
\rightarrow\infty,
\]
where $\nu_\alpha(x,\infty]=x^{-\alpha}$ for $x>0$.
The function $b(\cdot)$ is called a \textit{scaling function}.

If $\bX= (X_0,X_1,X_2,\ldots)$ is a (homogeneous) Markov chain and $K$
is a Markov transition kernel, we write $\bX\sim K$ to mean that the
dependence structure of $\bX$ is specified by $K$, i.e.,
\[
\P[X_{n+1} \in\cdot\given X_{n} = x ] = K(x,\cdot), \qquad
n=0,1,\ldots.
\]
Also, $\P_\mu[\bX\in\cdot]$ specifies the initial distribution $\P
[X_0 \in\cdot] = \mu$ (abbreviate $\P_x := \P_{\pp{x}}$), and
$\EP_\mu$ denotes expectation with respect to $\P_\mu$.

\section{Extremal component and tail chain approximation}
\label{secextrmch}

Let $\bX= (X_0,X_1,\ldots)$ be a Markov chain on $[0,\infty)$ with
transition kernel $K$.
The tail chain is a finite-dimensional approximation to the chain $\bX$
used to study the limit of $\{\Nn\}$ given by \eqref{eqppnn}.
Building on theory developed in \cite
{smith1992extremal,perfekt1997extreme,segers2007multivariate},
\cite{resnick2011asymptotics} presents the tail chain approximation in
terms of a related process known as the \textit{extremal component} of~$\bX$,
an approach we follow here.

\subsection{Tail chain approximation}

Suppose
the transition kernel
$K$ is in the domain of attraction of a distribution $G$ (denoted $K
\in D(G)$) which means \cite{resnick2011asymptotics},
\[
K(t,t\cdot) \wc G(\cdot) \qquad\mbox{on } [0,\infty] \qquad \mbox{as
} t\rightarrow\infty.
\]
Taking $\xi_1, \xi_2, \ldots$ i.i.d. random variables with common
distribution $G$, set $\xi(n) = \prod_{j=1}^n \xi_j$,\vspace*{2pt} $n\geq1$ with
$\xi(0)=1$ and write $\bxi=\{\xi(n), n\geq0\}$.
The \emph{tail chain associated with $G$} \cite
{perfekt1994extremal,smith1992extremal,resnick2011asymptotics}
is $\bT= (T_0,T_1,\ldots)$ given by
%
%
\begin{equation}
\label{eqtch} T_n = T_0 \xi_1 \cdots
\xi_n =T_0\xi(n),\qquad n\geq0.
\end{equation}
Thus $\bT$ is a multiplicative random walk
and $\{0\}$ is an absorbing barrier for $\bT$, accessible if $G(\{0\}
)> 0$.

An \emph{extremal boundary} for $\bX$ is a function $\eb(t)$
satisfying $0 \leq
\eb(t)\rightarrow0$, such that
%
%
\begin{equation}
\label{eqextrbdry} K\bigl(tu_t,t \bigl[0,\eb(t)\bigr]\bigr) \conv G
\bigl(\{0\}\bigr)\qquad \mbox{as } t\rightarrow\infty,
\end{equation}
for any non-negative function $u_t = u(t) \rightarrow u > 0$.
Such a function always exists if $K \in D(G)$ \cite
[Section~3.2]{resnick2011asymptotics}.
If $y(t)$ is an extremal boundary, any function $0\leq\tilde y(t)
\to0$ with $\tilde y(t) \geq y(t)$ for $t\geq t_0$ is also an extremal
boundary.
If $G(\{0\})= 0$, then $\eb(t)
\equiv0$ is a convenient choice.
Given an extremal boundary for $K$, the \textit{extremal component} of
$\bX$ is
the process $\bX$ prior to $\bX$ crossing below the scaled extremal
boundary and identically $\bzero$ afterwards. The first downcrossing
occurs at time
%
%
\begin{equation}
\label{eqn:tau} \tau(t) = \inf\bigl\{n \geq0 : X_n \leq t\eb(t)\bigr
\},
\end{equation}
and the extremal component is the process $\bX^{(t)} = (X_0^{(t)},
X_1^{(t)}, \ldots)$ defined by
\[
X_n^{(t)} = X_n \cdot\ind{n < \tau(t)},\qquad
n=0,1,\ldots.
\]
Starting from an extreme level $X_0=t$, the extremal boundary
separates extreme states from non-extreme states for
the scaled process $t^{-1}\bX$.

The \textit{tail chain approximation} is the following \cite
[Theorem~3.3]{resnick2011asymptotics}. For
$K\in D(G)$, $u>0$, $m\geq1$,
%
%
\begin{equation}
\label{eqecfddconv} \P_{tu} \bigl[ t^{-1} \bigl(
{X_1^{(t)}},\ldots,{X_m^{(t)}} \bigr)
\in\cdot \bigr] \wc\P_u \bigl[ (T_1,\ldots,
T_m ) \in\cdot \bigr] \qquad\mbox{on } [0,\infty]^m
\qquad\mbox{as } t\rightarrow\infty.
\end{equation}
So, the tail chain maps extreme states onto $(0,\infty)$ and
contracts non-extreme states to the point $\{0\}$.
Note $\tau(t) = \inf\{n \geq0 : X_n^{(t)} = 0 \}$.

If the finite-dimensional extremal behaviour of 
$\bX$
is completely accounted for by the extremal component, then the tail
chain approximation \eqref{eqecfddconv}
extends
from $\bX^{(t)}$ to $\bX$. When is this the case?
Say that $K$ satisfies the \emph{regularity condition} if
for any non-negative function $u_t = u(t) \rightarrow0$,
%
%
\begin{equation}
\label{eqregcond} K(tu_t,t\cdot) \wc\pp{0}(\cdot) \qquad \mbox{on }
[0,\infty] \qquad\mbox{as } t\rightarrow\infty.
\end{equation}
Equivalent forms of \eqref{eqregcond} exist in \cite
[Section~4]{resnick2011asymptotics}, and a relatively easy-to-check sufficient
condition is given in terms of update functions.
If either
(a) $y(t) \equiv0$ is an extremal boundary; or (b)
$K$ satisfies the regularity condition \eqref{eqregcond},
then for $u>0$, we strengthen \eqref{eqecfddconv} to \cite
[Theorem~4.1]{resnick2011asymptotics},
%
%
\begin{equation}
\label{eq:muscle} \P_{tu} \bigl[ t^{-1} (
{{X_1}},\ldots,{{X_m}} ) \in\cdot \bigr] \wc
\P_u \bigl[ (T_1,\ldots, T_m ) \in\cdot
\bigr] \qquad\mbox{on } [0,\infty]^m \qquad\mbox{as } t\rightarrow
\infty.
\end{equation}

\subsection{Finite-dimensional convergence}

The conditional approximation \eqref{eqecfddconv} requires that the
initial state become extreme. Combining \eqref{eqecfddconv} with a
heavy tailed initial distribution makes $\bX^{(t)}$ have an
unconditional distribution that is regularly
varying (in a sense to be discussed) with a limit measure
determined by the tail chain.
Depending on assumptions, convergences take place on
$\E_\sqsupset
:=(0,\infty]\times[0,\infty]^m$ or the larger
space $\E^*:=[0,\infty]^{m+1}\setminus\{\bzero\}$.

\begin{them}[({\cite[Proposition~5.1(b),
Theorem~5.1]{resnick2011asymptotics}})]\label{thmecjrv}
Let $\bX$ be a Markov chain on $[0,\infty)$ with $K\in D(G)$, and
suppose $X_0 \sim H$, with $1-H \in\RV_{-\alpha}$ with scaling
function $b(\cdot)$. On $\E_\sqsupset$ define the measure
%
%
\begin{equation}
\label{eqlimmeas} \mu(dx_0, d\boldsymbol x_m ) =
\nu_\alpha(dx_0)\P_{x_0} \bigl[ (T_1,
\ldots, T_m) \in d\boldsymbol x_m \bigr],
\end{equation}
and extend this to a measure $\mu^*$ on $\E^*$ by
defining
$\mu^*(\cdot\cap\E_\sqsupset) = \mu(\cdot)$ and
$\mu^*(\E^* \setminus\E_\sqsupset) = 0$.
For any $m\geq1$, the following convergences take place as
$t\rightarrow\infty$.
\begin{enumerate}[(b)]
\item[(a)]\label{parta}
In $\mplus((0,\infty]^m\times[0,\infty])$,
\[
P \bigl[(X_0,\ldots,X_m)/b(t) \in( \cdot)\cap
(0,\infty]^m\times[0,\infty] \bigr] \,\mathop{
\longrightarrow}^v\,\mu\bigl(( \cdot)\cap (0,
\infty]^m\times[0,\infty] \bigr),
\]
and in $\mplus(\E_\sqsupset)$
%
%
\begin{equation}
\label{eqecjrv} t\P \bigl[ \bigl({X_0^{(b(t))}}, \ldots,
{X_m^{(b(t))}} \bigr)/b(t) \in\cdot \bigr] \,\mathop{
\longrightarrow}^v\,\mu(\cdot).
\end{equation}
If either (\textup{i}) $G(\{0\})=0$; (\textup{ii}) $\eb(t)\equiv0$ is an extremal
boundary; or
(\textup{iii}) $K$ satisfies the regularity condition \eqref{eqregcond}, then
\eqref{eqecjrv} can be strengthened to
%
%
\begin{equation}
\label{eq:better} t \P \bigl[ ({X_0},\ldots,{X_{m}}
)/{b(t)} \in\cdot \bigr] \,\mathop{\longrightarrow}^v\, \mu(\cdot),
\qquad\mbox{in $\mplus(\E_\sqsupset)$.}
\end{equation}
\item[(b)]\label{partb}
In the bigger space $\E^*$, we have
%
%
\begin{equation}
\label{eqecjrvfull} t\P \bigl[ \bigl({X_0^{(b(t))}}, \ldots,
{X_m^{(b(t))}} \bigr)/{b(t)} \in\cdot \bigr] \,\mathop{
\longrightarrow}^v\,\mu^*(\cdot) \qquad\mbox{in }\mplus \bigl(\E^*
\bigr)
\end{equation}
if and only if
%
%
\begin{eqnarray}
\label{eqecmargrv} &&\EP\xi_1^\alpha< \infty\quad\mbox{and}
\quad t\P \bigl[X_j^{(b(t))} / b(t) \in\cdot \bigr] \,\mathop{
\longrightarrow}^v\, \bigl(\EP\xi_1^{ \alpha}
\bigr)^j \nua(\cdot)
\nonumber
\\[-8pt]
\\[-8pt]
&&\quad\mbox{in } \mplus(0,\infty],\qquad j=1, \ldots,m.
\nonumber
\end{eqnarray}
\end{enumerate}
\end{them}

Part (a) requires that the first observation is large and with
added conditions, part (b) removes this requirement.
Regardless of whether \eqref{eqecmargrv} holds, the limit is always a
lower bound on the tail weight of $X_j^{(t)}$, since
\[
\liminf_{t\rightarrow\infty} t\P \bigl[{X_j^{(b(t))}}/{b(t)} > x
\bigr] \geq\mu\bigl((0,\infty] \times[0,\infty]^{j-1} \times (x,\infty]
\times[0,\infty]^{m-j} \bigr) = \bigl(\EP\xi_1^{ \alpha}
\bigr)^j {x}^{-\alpha}
\]
by \eqref{eqecjrv} and Lemma~\ref{lemintmom} below.
Markov's inequality, \eqref{eqecjrv} and a moment condition:
%
%
\begin{equation}
\label{eqecmargrvmom} \exists\varepsilon> 0 \mbox{ such that }
\lim_{\delta\downarrow0} \limsup_{t\rightarrow\infty} t \EP \bigl[ \bigl(
{X_j^{(b(t))}}/{b(t)} \bigr)^\varepsilon\ind{X_0
\leq\delta b(t) } \bigr] = 0,\qquad j=1,\ldots,m,
\end{equation}
imply \eqref{eqecmargrv}. See
\cite{maulik2002asymptotic}.

Here is a formula that helps evaluate $\mu$ in \eqref{eqlimmeas} on
certain sets.

\begin{lem} \label{lemintmom}
For a random variable $Y$, define the measure $\nu(dx,dy) = \nua(dx)
\P[x Y \in dy]$ on $[0,\infty]^2 \setminus\{(0,0)\}$.
We compute
%
%
\begin{equation}
\label{eqintmom} \nu\bigl( [0,x] \times(y,\infty] \bigr) = {y}^{-\alpha} \EP \bigl[
Y^\alpha\ind{Y > yx^{-1}} \bigr] - {x}^{-\alpha} \P \bigl[
Y > yx^{-1} \bigr].
\end{equation}
In particular, $\nu([0,\infty] \times(y,\infty]) = {y}^{-\alpha}
\EP
Y^\alpha$.
\end{lem}

\begin{pf}
We obtain
\begin{eqnarray*}
\nu\bigl([0,x]\times(y,\infty] \bigr) &=& \int_{[0,x]} \nua(ds) \P
\bigl[Y > y s^{-1} \bigr] = \int_{[{x}^{-\alpha},\infty]} ds\, \P \bigl[Y
> y s^{1/\alpha} \bigr]
\\
&=& \int_{[{x}^{-\alpha},\infty]} ds\, \P \bigl[{y}^{-\alpha}
Y^\alpha> s \bigr]
\end{eqnarray*}
by change of variables. Applying Fubini's theorem, this becomes
\[
\int_{({x}^{-\alpha},\infty]} \bigl(s-{x}^{-\alpha} \bigr) \P
\bigl[{y}^{-\alpha} Y^\alpha\in ds \bigr] = {y}^{-\alpha} \EP
\bigl[ Y^\alpha\ind{Y > yx^{-1}} \bigr] - {x}^{-\alpha} \P
\bigl[ Y > yx^{-1} \bigr].
\]
Letting $x\rightarrow\infty$, this quantity converges to
${y}^{-\alpha}
\EP
Y^\alpha$ by monotone convergence.
\end{pf}


\subsection{Maximum of the extremal component}

We give conditions on the extremal component which enable an
informative
point process limit result by controlling the positive portion of the
extremal component, the random vector of random length
$ \{X_j^{(t)}:j=0,\ldots,\tat-1\} = \{{X_j} : j=0,\ldots,\tat-1\}$.
The conditions imply restrictions on the behaviour of the tail chain
$\bT$.

We study a positive recurrent chain $\bX$ by splitting it into
regenerative cycles and analyzing its extremal properties via the
extremal components of the cycles.
For regenerative processes, Asmussen
\cite{asmussen1998extreme} and Rootz\'en \cite{rootzen1988maxima}
point out the connection between point process convergence and the asymptotic
distribution of cycle maxima. Informed by this approach, we
consider when the distribution of the maximum
over the extremal component has a regular variation property.
The limit measure of this regular variation can be used to compute an extremal
index for $\bX$ \cite{rootzen1988maxima}.

Here is a condition that controls the persistence of non-zero values of
the extremal component:
%
%
\begin{equation}
\label{condeclrd} \lim_{m\rightarrow\infty} \limsup_{t\rightarrow
\infty} \P \Bigl[
\sup_{j
\geq m} {X_j^{(b(t))}}/{b(t)} > a \big|
{X_0} > \delta{b(t)} \Bigr] = 0 \qquad\mbox{for all } a, \delta> 0.
\end{equation}
Note\vspace*{1pt} $ \sup_{j\geq m} X_j^{(b(t))} = (\sup_{m \leq j < \tau
(b(t))} X_j
) \ind{\tau(b(t)) > m}$.
Compare this condition with \cite[Condition 4.1]{basrak2009regularly}
and \cite[Equation (3.1)]{perfekt1994extremal}, which are formulated
in terms
of block sizes.
Condition \eqref{condeclrd} is a tightness condition that complements
the finite-dimensional convergences \eqref{eqecjrv}. Section~\ref
{secmch} gives simpler sufficient conditions
depending on whether $G(\{0\})> 0$ or $G(\{0\})= 0$.
Condition \eqref{condeclrd} requires
that the chain drift back to the non-extreme states after visiting an
extreme state and makes non-extreme states
recurrent and the tail chain
transient.

\begin{prop} \label{propcondeclrd}
Let $\bX\sim K\in D(G)$ be Markov on $[0,\infty)$ with
initial distribution $H$ satisfying $1-H \in\RV_{-\alpha}$ with
scaling function $b(t)$, so that
\eqref{eqecjrv} holds. If $\bX$ satisfies Condition
\eqref{condeclrd}, then
%
%
\begin{equation}
\label{eqcondeclrdtch} \lim_{m\rightarrow\infty} \P \Bigl[ \sup_{j \geq m} \xi(j)
> a \Bigr] = 0,\qquad a>0,
\end{equation}
and $\xi(n) \to0$ as $n \to\infty$ in probability and almost surely
and therefore,
%
%
\begin{equation}
\label{eqtchtransient} \P[T_n \rightarrow0] = 1.
\end{equation}
\end{prop}

So the tail chain is transient and the
additive random walk $\{ \log T_n \}_{n \geq0}$ satisfies
$\log T_m \to-\infty$.
The tail chain $\bT$ and $\bX$ live on the same state space
$[0,\infty)$ but for $\bT$, $\{0\}$ is a special boundary
state which represents the collection of non-extreme states of $\bX$
under the tail chain approximation.

\begin{pf*}{Proof of Proposition~\ref{propcondeclrd}}
Observe from \eqref{eqecjrv}, as $t\to\infty$,
\begin{eqnarray*}
t\P \Bigl[ X_0 > b(t), \sup_{m\leq j\leq r} {X_j^{(b(t))}}
> {b(t)} \Bigr] &\conv& \int_{(1,\infty]} \nua(dx) \P_x
\Bigl[\sup_{m\leq j\leq r} T_j > 1 \Bigr]
\\
&=& \int_{(1,\infty]} \nua(dx) \P \Bigl[\sup_{m\leq j\leq r} \xi (j)
> x^{-1} \Bigr]. 
\end{eqnarray*}
Therefore, by monotonicity and then monotone convergence,
\begin{eqnarray*}
\limsup_{t\rightarrow\infty} t\P \Bigl[ X_0 > b(t), \sup_{j \geq m}
{X_j^{(b(t))}} >{b(t)} \Bigr] 
&\geq&\lim_{r \rightarrow\infty} \int
_{(1,\infty]} \nua(dx) \P \Bigl[ \sup_{m \leq j \leq r} \xi(j) >
x^{-1} \Bigr]
\\
&= &\int_{(1,\infty]} \nua(dx) \P \Bigl[ \sup_{j \geq m} \xi(j)
> x^{-1} \Bigr]
\\
&=:& \int_{(1,\infty]} \nua(dx) f_m(x).
\end{eqnarray*}
Condition\vspace*{1pt} \eqref{condeclrd} implies that $\int_{(1,\infty]} \nua
(dx) f_m(x) \rightarrow0$ as $m\rightarrow\infty$.
We claim that $f_m(x) \rightarrow0$ for any $x>0$.
Suppose instead that $\inf_m f_m(x_0) \geq c > 0$ for some $x_0$.
Since the $f_m$ are all increasing in $x$, we have $\inf_m f_m(x) \geq
c$ for $x\geq x_0$.
But this implies that
\[
\liminf_{m\rightarrow\infty} \int_{(1,\infty]} \nua(dx) f_m(x)
\geq\liminf_{m\rightarrow\infty} \int_{(1 \smax x_0,\infty]} \nua(dx)
f_m(x) \geq c \nua(1 \smax x_0,\infty] > 0
\]
by Fatou's Lemma, contradicting Condition~\ref{condeclrd}. Therefore,
$\P[ \sup_{j \geq m} \xi(j) > x^{-1} ] \rightarrow0$ as
$m\rightarrow
\infty$ for all $x>0$, establishing \eqref{eqcondeclrdtch}.
\end{pf*}

Condition \eqref{condeclrd} assumes the first observation exceeds
$\delta b(t)$ which is in the spirit of \eqref{eqecjrv}.
For translating the stronger convergence of unconditional
distributions \eqref{eqecjrvfull} in the bigger space
$\E^*:=[0,\infty]^{m+1}\setminus\{\bzero\}$ to point process
convergence, we will require an
additional assumption:
%
%
\begin{eqnarray}
\label{condeclrdmrv} &&\exists m_0 \geq1 \mbox{ such that}
\nonumber
\\[-8pt]
\\[-8pt]
&&\quad\lim_{\delta\downarrow0} \limsup_{t\rightarrow\infty} t\P \Bigl[ {X_0}/{b(t)}
\leq\delta, \sup_{j\geq m_0} {X_j^{(b(t))}}/{b(t)} > a \Bigr] = 0 \qquad\mbox{for all } a > 0.
\nonumber
\end{eqnarray}
Analogously to \eqref{eqecmargrvmom}, by Markov's inequality a moment
condition
is sufficient for Condition \eqref{condeclrdmrv}:
\[
\lim_{\delta\downarrow0} \limsup_{t\rightarrow\infty} t \EP \Bigl[ \Bigl(
\sup_{j\geq m_0} {X_j^{(b(t))}}/{b(t)} \Bigr)^\varepsilon
\ind{X_0 \leq\delta b(t)} \Bigr] = 0, \qquad\mbox{for some $
\varepsilon> 0$}.
\]
Condition \eqref{condeclrdmrv} implies a
uniform bound on the $\alpha$th moment of the tail chain states.

\begin{prop}\label{pro2.2}
Let $\bX\sim K\in D(G)$ be a Markov chain on $[0,\infty)$ with
initial distribution $H$ satisfying $1-H \in\RV_{-\alpha}$, whose
extremal component satisfies \eqref{eqecjrvfull} on $\mplus(\E^*)$.
If $\bX$ satisfies Condition
\eqref{condeclrdmrv} then,
%
%
\begin{equation}
\label{eqcondeclrdtchmom} \EP \Bigl( \sup_{j \geq1} \xi(j)^\alpha
\Bigr) < \infty.
\end{equation}
\end{prop}

\begin{rmk*}
Under \eqref{eqcondeclrdtchmom}, we necessarily have $\EP\xi_1^\alpha
\leq1$ since
\[
\sup_{j \geq1} \bigl(\EP\xi_1^\alpha
\bigr)^j = \sup_{j \geq1} \EP\xi(j)^\alpha\leq\EP \Bigl(
\sup_{j \geq1} \xi(j)^\alpha \Bigr) < \infty.
\]
Recalling \eqref{eqecmargrv}, the marginal tails of
the extremal component cannot be heavier than the tail of~$H$.
\end{rmk*}

\begin{pf*}{Proof of Proposition~\ref{pro2.2}}
Under \eqref{eqecjrvfull},
\[
t\P \Bigl[ \sup_{m_0 \leq j \leq r} {X_j^{(b(t))}} > b(t) \Bigr]
\conv\int_{(0,\infty]} \nua(dx) \P \Bigl[ \sup_{m_0 \leq j \leq
r} \xi(j)
> x^{-1} \Bigr] 
= \EP \Bigl(\sup_{m_0 \leq j \leq r} \xi(j)^\alpha \Bigr)
\]
by Lemma~\ref{lemintmom}. Therefore,
\[
\limsup_{t\rightarrow\infty} t\P \Bigl[ \sup_{j\geq m_0} {X_j^{(b(t))}}
> {b(t)} \Bigr] \geq\lim_{r\rightarrow\infty} \EP \Bigl(\sup_{m_0 \leq j \leq r} \xi
(j)^\alpha \Bigr) = \EP \Bigl(\sup_{j \geq m_0} \xi(j)^\alpha
\Bigr).
\]
Furthermore, by Condition \eqref{condeclrdmrv}, for some $\delta>0$,
\begin{eqnarray*}
&&\limsup_{t\rightarrow\infty} t\P \Bigl[ \sup_{j\geq m_0} {X_j^{(b(t))}}
> {b(t)} \Bigr]
\\
&&\quad\leq\limsup_{t\rightarrow\infty} t\P \Bigl[ X_0 \leq \delta b(t),
\sup_{j\geq m_0} {X_j^{(b(t))}} > {b(t)} \Bigr] +
\limsup_{t\rightarrow
\infty} t\P \bigl[ X_0 > \delta{b(t)} \bigr] < \infty
\end{eqnarray*}
showing that $\EP(\sup_{j\geq m_0} \xi(j)^\alpha) <
\infty$
which is enough for \eqref{eqcondeclrdtchmom}.
\end{pf*}

Under both Conditions \eqref{condeclrd} and \eqref{condeclrdmrv},
we derive the tail behaviour of the maximum of the extremal component
of $\bX$.

\begin{prop} \label{proprvecmax}
Let $\bX\sim K\in D(G)$ be a Markov chain on $[0,\infty)$ with
initial distribution $H$ satisfying $1-H \in\RV_{-\alpha}$, whose
extremal component satisfies both \eqref{eqecjrvfull} on $\mplus(\E^*)$.
If $\bX$ satisfies Conditions
\eqref{condeclrd} and \eqref{condeclrdmrv}, then
%
%
\begin{equation}
\label{eqrvcyc} t\P \Bigl[\sup_{0 \leq j < \tau(b(t))} {X_j}/{b(t)} \in
\cdot \Bigr] \,\mathop{\longrightarrow}^v\,c \cdot\nua(\cdot) \qquad
\mbox{in } \mplus(0,\infty],
\end{equation}
where
$ c = \P[ \sup_{j \geq1} \xi(j) \leq1 ] + \EP[\sup_{j
\geq1} \xi(j)^\alpha\ind{\sup_{j \geq1} \xi(j) > 1} ]$.
\end{prop}

\begin{pf}
For $x>0$, we have $[\sup_{j<\tau(b(t))}X_j /b(t) >x] =[\sup_{j\geq
0} X_j^{(b(t))}/b(t)>x]$.
For $m\geq1$, we have on the one hand, by \eqref{eqecjrvfull},
\begin{eqnarray*}
\liminf_{t\rightarrow\infty} t\P \Bigl[\sup_{0 \leq j < \tau(b(t))} {X_j}/{b(t)} >
x \Bigr] &\geq& \lim_{t\rightarrow\infty} t\P \Bigl[\sup_{0 \leq
j < m}
{X_j^{(b(t))}}/{b(t)} > x \Bigr]
\\
&=& {x}^{-\alpha} + \int_{[0,x)} \nua(du)
\P_u \Bigl[ \sup_{1 \leq
j <
m} T_j > x \Bigr],
\end{eqnarray*}
from which, letting $m\rightarrow\infty$,
%
%
\begin{equation}
\label{eqpfeclb} \liminf_{t\rightarrow\infty} t\P \Bigl[\sup_{0
\leq j < \tau(b(t))}
{X_j}/{b(t)} > x \Bigr] \geq{x}^{-\alpha} + \int
_{[0,x)} \nua(du) \P_u \Bigl[ \sup_{j
\geq1}
T_j > x \Bigr].
\end{equation}
On the other hand, for $\delta>0$ we have
\[
t\P \biggl[\sup_{j \geq m} \frac{X_j^{(b(t))}}{b(t)} > x \biggr] \leq t\P \biggl[
\frac{{X_0}}{b(t)} > \delta, \sup_{j \geq m} \frac
{X_j^{(b(t))}}{b(t)} > x \biggr] + t
\P \biggl[\frac{{X_0}}{b(t)} \leq\delta, \sup_{j \geq m} \frac
{X_j^{(b(t))}}{b(t)} > x
\biggr].
\]
Given $\varepsilon> 0$, by
Condition \eqref{condeclrdmrv}, we may choose $\delta$ small enough that
\[
\limsup_{t\rightarrow\infty} t\P \Bigl[{{X_0}} \leq\delta{b(t)},
\sup_{j \geq m_0} {X_j^{(b(t))}} > {b(t)}x \Bigr] <
\varepsilon/2,
\]
where $m_0$ is from Condition \eqref{condeclrdmrv}.
Condition
\eqref{condeclrd} permits the choice $m_1 \geq m_0$ so large that
\[
\limsup_{t\rightarrow\infty} t\P \Bigl[{{X_0}} > \delta{b(t)},
\sup_{j
\geq m_1} {X_j^{(b(t))}} > {b(t)}x \Bigr] <
\varepsilon/2.
\]
Therefore,
$ \limsup_{t\rightarrow\infty} t\P[\sup_{j \geq m}
{X_j^{(b(t))}}/{b(t)} > x ] < \varepsilon$ for $m \geq m_1$,
and so
\begin{eqnarray*}
\limsup_{t\rightarrow\infty} t\P \Bigl[\sup_{0 \leq j < \tau(b(t))} {X_j}/{b(t)} >
x \Bigr] &< &\lim_{m\rightarrow\infty} \limsup_{t\rightarrow
\infty} t\P \Bigl[
\sup_{0
\leq j <m} {X_j^{(b(t))}}/{b(t)} > x \Bigr] +
\varepsilon
\\
&=& {x}^{-\alpha} + \int_{[0,x)} \nua(du)
\P_u \Bigl[ \sup_{j \geq1} T_j > x \Bigr] +
\varepsilon.
\end{eqnarray*}
Combine this with \eqref{eqpfeclb}, and apply formula
\eqref{eqintmom} for $\nu([0,x]\times(x,\infty])$ to complete the proof.
\end{pf}

\section{Point process convergence for Markov chains}
\label{secmch}

We now derive the limit of the exceedance point process $\Nn$ defined
in \eqref{eqppnn}, where $\bX= (X_0,X_1,\ldots)$ is a Markov chain
on $[0,\infty)$ with transition kernel $K \in D(G)$.
We write
%
%
\begin{equation}
\label{eqn:Nn} \Nn= \sum_{j=0}^\infty
\pp{({j}/{n}, {X_j}/{b_n})},
\end{equation}
using the notation $\pp{x}$ to denote the measure assigning unit mass
at the point $x$ and $\Nn$ is a random element of
$\Mp([0,\infty)\times(0,\infty])$, the space of Radon point measures
on $[0,\infty)\times(0,\infty]$, endowed with the topology of vague
convergence \cite{resnick1987extreme,kallenberg1976random,neveu1977processus}.

If $\bX$ is
positive recurrent, it is a regenerative process (\cite[Section
VII.3]{asmussen2003applied}, \cite{meyn2009markov}) so the sample
path of $\bX$ splits into identically distributed cycles
between visits to certain set. The extremal properties of $\bX$
are determined by extremal behaviour of the individual cycles. This
approach has
been developed for Markov chains by Rootz\'en \cite
{rootzen1988maxima}, as
well as for queues by Asmussen
\cite{asmussen1998extreme}. Our approach
introduces the tail chain approximation to describe the extremal
behaviour of the regenerative cycles using their extremal component.

\subsection{Cycle decomposition}

Consider the case where $\bX$ has a positive recurrent atom $A$. For
positive recurrent chains, atoms can be constructed by several methods
if no natural atom exists. See, e.g.,\vadjust{\goodbreak} \cite[Chapter~6]{durrett:1991} or
\cite[Chapter I.5]{meyn:tweedie:1993}.
An atom is a set such that for a probability distribution $H$ on
$[0,\infty)$,
%
%
\begin{equation}
\label{eqatom} K(y,\cdot) = H(\cdot) \qquad\mbox{for all } y \in A \quad
\mbox{and}\quad
\P_y[\ta< \infty] = 1 \qquad\mbox{for } y
\geq0,
\end{equation}
and $\ta= \inf\{n\geq0 : X_n \in A \}$ is the first hitting time of $A$.
Positive recurrence means that
%
%
\begin{equation}
\label{eqposrec} \EP_H \ta< \infty,
\end{equation}
where $\EP_H$ denotes expectation with respect to $H$ considered as
the initial distribution of~$X_0$.

Under \eqref{eqatom}, the sample path of $\bX$ splits into i.i.d.
cycles between visits to $A$, as follows.
Define the times $\{S_k\}$, $\{\tau^A_k\}$ recursively according to
%
%
\begin{eqnarray}
\label{eqrenewaltimes} &&\tau^A_0= \ta,
\hspace*{108.5pt} S_0 = \tau^A_0+ 1;
\\
&&\tau^A_k= \inf\{n \geq0 : X_{S_{k-1} + n} \in A \},
\qquad S_k = S_{k-1} + \tau^A_k+
1, k \geq1.
\nonumber
\end{eqnarray}
Thus, the sequence $0 \leq S_0 - 1 < S_1 - 1 < S_2 - 1 < \cdots$
gives the indices when $\bX$ is in $A$, and $X_{S_k} \sim H$ for
$k\geq0$. The values $\tau^A_k\geq0$ are the number of steps
$\bX$ takes outside of $A$ between visits to $A$.
The cycles end by visits to $A$; cycles are the random elements
\[
C_0 = 
(X_0, X_1, \ldots,
\mathop{\mathop{X_{\tau^A_0}}_{\in}}_{A}) 
\quad\mbox{and}\quad C_k = 
(\mathop{\mathop{X_{S_{k-1}}}_{\wr}}_{H},
\ldots, \mathop{\mathop{X_{S_{k-1} + \tau^A_k}}_{\in}}_{A})
,\qquad k\geq1
\]
of the space of finite sequences $\mathcal{S} = \bigcup_{m=1}^\infty
\R^m$.
The strong Markov property implies
$C_0,C_1,\ldots$ are independent, and $C_1, C_2,
\ldots$ are identically distributed.
In particular, for $k\geq1$,
\[
\P \bigl[ \bigl\{C_k ; \tau^A_k\bigr\} \in
\cdot \bigr] = \P \bigl[ \bigl\{ (X_{S_{k-1}}, \ldots, X_{S_{k-1} + \tau^A_k}) ;
\tau^A_k \bigr\} \in \cdot \bigr] = \P_H
\bigl[ \bigl\{(X_0,\ldots, X_{\ta}) ; \ta \bigr\} \in\cdot
\bigr].
\]
Furthermore, $0 < S_0 < S_1 < S_2 < \cdots$ is a renewal process, with
%
%
\begin{equation}
\label{eqdefq} q = \EP(S_1-S_0) = \EP_H
\tau_A + 1 < \infty
\end{equation}
by \eqref{eqposrec}.
Applying the cycle decomposition, we may now write \eqref{eqn:Nn} as
%
%
\begin{equation}
\label{eqnndecomp} \Nn= \sum_{j=0}^{\infty}
\pp{({j}/{n}, {X_j}/{b_n} )} 
= \sum_{j=0}^{S_0-1} \pp{({j}/{n},
{X_j}/{b_n} )} 
+ \sum
_{k=0}^\infty\sum_{j=0}^{\tau^A_{k+1}}
\pp{ ({(S_k+j)}/{n}, {X_{S_k+j}}/{b_n}
)} 
= \chi_n^0+ \chi^*_n.
\end{equation}
%
As a family
of random elements in $\Mp
([0,\infty)\times(0,\infty])$, $\{\chi_n^0\}$ is asymptotically negligible.

\begin{lem} \label{propppneg}
Assuming \eqref{eqposrec} and $b_n \to\infty$,
$\chi_n^0\Rightarrow0$, the null measure, in
$\Mp
([0,\infty)\times(0,\infty])$.
\end{lem}

\begin{pf}
Let $f\in\mathcal{C}^{+}_K([0,\infty) \times(0,\infty])$ with
support in
$[0,R]\times[M,\infty]$ for integers $R,M$.
It is sufficient to verify that $\P[\chi_n^0(f) > \gamma]
\rightarrow0$,
for any $\gamma>0$.
We have as $n\rightarrow\infty$,
\begingroup
\abovedisplayskip=6.5pt
\belowdisplayskip=6.5pt
\begin{eqnarray*}
\P \bigl[\chi_n^0(f) > \gamma \bigr] &=& \P \Biggl[ \sum
_{j=0}^{S_0-1} f \biggl(\frac{j}{n},
\frac{X_j}{b_n} \biggr) > \gamma \Biggr]
\\[-2pt]
&\leq&\sum_{m=0}^r \P \Biggl[ \sum
_{j=0}^{m} f \biggl( \frac{j}{n},
\frac{X_j}{b_n} \biggr) > \gamma, \ta= m \Biggr] + \P[\ta> r ],
\end{eqnarray*}
and
\begin{eqnarray*}
\sum_{m=0}^r \P \Biggl[ \sum
_{j=0}^{m} f \biggl(\frac{j}{n},
\frac
{X_j}{b_n} \biggr) > \gamma, \ta= m \Biggr] &\leq&(r+1) \P \Biggl[ \sum
_{j=0}^{r} f \biggl(\frac{j}{n},
\frac
{X_j}{b_n} \biggr) > \gamma \Biggr]
\\[-2pt]
&\leq&(r+1) \P \Bigl[ \sup_{0 \leq j \leq r} X_j \geq b_n
M \Bigr] \conv0.
\end{eqnarray*}
Choosing $r$ to make $\P[\ta> r]$ arbitrarily small,
the result follows.
\end{pf}

\subsection{Point process convergence}

Lemma~\ref{propppneg} and Slutsky's theorem means that the asymptotic
behavior of
$N_n$ and $\chi_n^*$ are the same.
We obtain a weak limit for $\chi_n^*$
using the tail chain approximation discussed in Section~\ref
{secextrmch}, provided that a cycle's extremal behaviour is
adequately described by its extremal component.

As usual, assume $K \in D(G)$, $1-H\in\RV_{-\alpha}$
and suppose $y(t)$ is an extremal boundary for $\bX$.
We require a mild assumption that the atom $A$ be a bounded subset of
the state space $[0,\infty)$,
%
%
\begin{equation}
\label{eqatombdd} \sup A < \infty,
\end{equation}
as would usually be the case. Fix $k\geq1$. The number of steps needed
by the scaled
process
in the $k$-th cycle to cross below the extremal boundary is
\[
\tau_{k}(t) = \inf \bigl\{ n\geq0 : X_{S_{k-1} + n} \leq t y(t)
\bigr\}.
\]
The extremal component of the $k$-th cycle is
$C_k(t) := \{X_{S_{k-1}+j} : j=0,\ldots,\tkt-1 \}$.

Without loss of generality, we suppose the extremal component of a
cycle is a subset of the complete cycle. To see this, observe from the
definition of $\tau(t) $ and $\tau_A$, without
loss of generality,
%
%
\begin{equation}
\label{eqcycletimes} \P\bigl[\tat\leq\ta, \forall t>0 \bigr] = 1.
\end{equation}
\endgroup
Indeed, \eqref{eqatombdd} implies that $A \subset[0,c]$ some $c$.
Define ${\tau}_c = \inf\{n\geq0 : X_n \leq c\}$ and
$\P[{\tau}_c \leq\ta]=1$; we claim further that we may suppose $\P
[\tat\leq{\tau}_c, \forall t>0]=1$.
If $y(t) \geq c/t$ for all $t>0$, then this follows directly.
Otherwise, verify that $\tilde{y}(t) = y(t) \smax c/t$ is also an
extremal boundary for $K$ (see the remarks after \eqref{eqextrbdry}),
and the corresponding downcrossing time
satisfies $\P[\tilde{\tau}(t) \leq{\tau}_c, \forall t>0]=1$.\vadjust{\goodbreak}

Therefore, for $k\geq1$,
%
%
\begin{equation}
\label{eqcycdistrib} 
\hspace*{-5pt}\P \bigl[ \bigl\{ (X_{S_{k-1}}, \ldots,
X_{S_{k-1} + \tkt-1}) ; \tkt, \tau^A_k \bigr\} \in\cdot
\bigr] = \P_H \bigl[ \bigl\{(X_0,\ldots,
X_{\tat-1}) ; \tat, \ta \bigr\} \in \cdot \bigr]
\end{equation}
and $\{(C_k(t); \tkt, \tau^A_k) : k \geq1 \}$ are independent, since
each is a function of $\{C_k; \tau^A_k\}$.
These facts suggest we approximate $\chi^*_n$ by a point process whose
observations consist of the extremal components of iid
copies of the chain $\bX$ started from $X_0 \sim H$. This
approximation is facilitated by additional notation.
Let $\{\bX,\bX_k=(X_{kj}, j\geq0):k\geq0\}$ be i.i.d. copies of the
Markov chain $\bX\sim K$ with respect to $\P_H (\cdot)$; that is the initial
distribution of each chain is
fixed to be
$H$.
Define
\[
\tkoit= \inf \bigl\{ j\geq0 : X_{kj} \leq t y(t) \bigr\},\qquad k=0,1,
\ldots,
\]
and for $k\geq0$, form the extremal component $\bX_k^{(t)} = \{X_{kj}
\cdot\ind{j < \tkoit}, j\geq0 \} $ of the $k$th chain.
Thus with respect to $\P_H(\cdot)$, $ (\bX^{(t)}_k, \tkoit)
\,\mathop{=}\limits^{d}\,
(\bX^{(t)}, \tau(t)
)$ for $k\geq0$, with the tilde differentiating the times
$\tkit$ defined on the $k$th process $\bX_{k}$ from the cycle times
$\tkt$ defined on $\bX$. Recall $\tau(t) $ is also defined on $\bX$.

Next, generate an i.i.d. family of tail chains by letting $\{\bxi,
\bxi_k = (\xi_k(n), n\geq0):\break k\geq0 \}$ be i.i.d. copies
of the process $\bxi= (\xi(n), n\geq0)$, recalling the notation
around \eqref{eqtch}.
Additionally, put $\tau^*_{k+1}= \inf\{ j\geq0 : \xikj= 0\}$, the first
time the $k$th tail chain hits $0$. Use
the convention $\inf\varnothing= \infty$; for example, $\tau^*_{k+1}=
\infty$ a.s., $k\geq0$, if $G(\{0\})= 0$.
Finally, let
\[
\limpp= \sum_{k=0}^\infty
\pp{(t_k, i_k)} \sim\prm(\leb\times\nua),
\]
be a Poisson random measure on $\Mp([0,\infty)\times(0,\infty])$,
independent of the $\{\bxi_k\}$, with mean measure
a product of Lebesgue measure
on the time axis $[0,\infty)$ and Pareto measure $\nua$ (given by
$\nua(x,\infty] = {x}^{-\alpha}$) on the observation axis $(0,\infty]$.
Recall $\alpha$ is the tail index of $\bar H$.

The point process consisting of the observations $\bX_k^{(b_n)}$, spaced
in time according to the renewal times $\{S_k\}$, converges to a
cluster Poisson process which is basically $\zeta$ with a time
scaling and compounded in the second coordinate
according to the i.i.d. tail chains $\{\bxi_k\}$.
This result is basic to analyzing the asymptotic behavior of $N_n$ in
\eqref{eqn:Nn}.

\begin{prop} \label{thmppconv}
Let $\bX$ be a Markov chain on $[0,\infty)$ with transition kernel
$K\in D(G)$, and initial distribution $ H$, such that
$tH(b(t)\cdot) \,\mathop{\rightarrow}\limits^v\,
\nu_\alpha(\cdot)$ in $\mplus(0,\infty]$, where $b(t)\rightarrow
\infty$
and $b_n=b(n)$. The renewal process
$\{S_k\}$ is defined in
\eqref{eqrenewaltimes}, with mean interarrival time $q$ given by
\eqref{eqdefq}. With the notation introduced in the previous
paragraphs, we have the following with respect to $\P_H$.
\begin{enumerate}[(b)]
\item[(a)] If $\bX$ satisfies Condition \eqref{condeclrd}, then given
$\delta>0$,
in $\Mp([0,\infty)\times(0,\infty])$, as $n\rightarrow\infty$,
%
%
\begin{equation}
\label{weenie1} \ppnmain:= \sum_{k=0}^\infty
\sum_{j=0}^{\tk-1} \pp{(
\frac
{S_k+j}{n}, \frac{\Xkj}{b_n} )} \ind{\Xko\geq\delta b_n
} \wc\sum_{k=0}^\infty\sum
_{j=0}^{\tau^*_{k+1}-1} \pp{(qt_k, i_k
\xikj)} \ind{i_k \geq\delta} =: \ppmain.\vadjust{\goodbreak}
\end{equation}

\item[(b)] Suppose ${\bX}$ satisfies
\eqref{eqecjrvfull} as well as both
Conditions \eqref{condeclrd} and \eqref{condeclrdmrv}. Then
in $\Mp([0,\infty)\times(0,\infty])$, as $n\rightarrow\infty$,
%
%
\begin{equation}
\label{weenie2} \eta^*_n:= \sum_{k=0}^\infty
\sum_{j=0}^{\tk-1} \pp{(
\frac
{S_k+j}{n}, \frac{\Xkj}{b_n} )} 
\wc\sum
_{k=0}^\infty\sum_{j=0}^{\tau^*_{k+1}-1}
\pp{(qt_k, i_k\xikj)} =: \eta^*.
\end{equation}
\end{enumerate}
\end{prop}

Section~\ref{secppconv} (p. \pageref{secppconv}) contains the proof.
Paralleling the discussion in Section~\ref{secextrmch}, we
have two results depending on the strength of the
conditions. The weaker assumptions of part (a)
yield a
result that selects cycles starting from an
exceedance. Part (b) does not have to do such cycle selection.

The points of the limit process are arranged in stacks above
common time points $qt_k$. The heights of the points in each stack are
specified by an independent run of the tail chain starting from
$i_k$. If $G(\{0\})> 0$, then the $\tau_k^*$ are i.i.d. Geometric random
variables with parameter $G(\{0\})$, so all stacks have finite
length. If $G(\{0\})=0$, then $\P[\tau_k^* = \infty]=1$ for each
$k$. In
this case, Condition \eqref{condeclrd} is necessary to ensure that
$\eta^*$ is Radon, by forcing the tail chain to drift towards 0 as in
\eqref{eqtchtransient}. The process $\ppmain$ retains only those
stacks of $\eta^*$ whose initial value exceeds the threshold
$\delta$. Because there are an infinite number of $i_k$ in any
neighbourhood of 0, dispensing with the restriction in $\delta$
requires that not too many of the $\xikj$ are large. This translates
to the condition $\EP\xi_1^\alpha\leq1$, provided by Condition~\eqref{condeclrdmrv}.

To analyze $N_n$ in \eqref{eqn:Nn},
we approximate $\chi^*_n$ in \eqref{eqnndecomp} by $\eta^*_n$ in
\eqref
{weenie2}, provided the extremal
component adequately describes extremal behaviour within each
cycle.
If the extremal boundary is not identically zero, behavior
between the end of the extremal component and the end of
the cycle is not be captured by the tail chain and we require
that such observations do not significantly
influence extremal properties. To guarantee a result analogous to
Part (a)
above, we require,
%
%
\begin{equation}
\label{condcyc} \lim_{t\rightarrow\infty} \P \Bigl[ \sup_{\tau(b(t))<
j < \ta}
{{X_j}}/{b(t)} > a \big| X_0 > \delta b(t) \Bigr] = 0 \qquad
\mbox{for all } a, \delta> 0,
\end{equation}
and for a result analogous to Part (b) above, we require,
%
%
\begin{equation}
\label{condcycmrv} \lim_{t\rightarrow\infty} t\P \Bigl[ \sup_{\tau(b(t))< j < \ta}
{{X_j}}/{b(t)} > a \Bigr] = 0 \qquad\mbox{for all } a > 0.
\end{equation}
With these conditions, the point process $\Nn$ converges to the limit
$\eta^*$, and the distribution of the cycle maximum behaves as if it
has a regularly
varying tail.

\begin{them} \label{thmppconvmc}
Let $\bX$ be a Markov chain on $[0,\infty)$ with transition kernel
$K\in D(G)$.
Suppose that $K$ has a positive recurrent bounded atom in the sense of
\eqref{eqatom}, \eqref{eqposrec}, and~\eqref{eqatombdd}. Define the
renewal process $\{S_k\}$ with mean interarrival time $q$ as in
\eqref{eqrenewaltimes} and \eqref{eqdefq} and assume further that
$tH(b(t)\cdot) \,\mathop{\rightarrow}\limits^v\,\nu_\alpha(\cdot
)$ in
$\mplus(0,\infty]$, where $b(t)\rightarrow\infty$. With respect to
$\P_H$,
the following hold.
\begin{enumerate}[(b)]
\item[(a)]
If $\bX$ satisfies Conditions \eqref{condeclrd} and \eqref{condcyc},
then given $\delta>0$,
%
%
\begin{equation}
\label{eqppconvxdfull} 
\Nnd:= \sum_{0\leq j<S_0}
\pp{(\frac{j}{n}, \frac{X_j}{b_n} )} +\sum
_{k=1}^{\infty} \sum_{S_{k-1} \leq j < S_k}
\ind{X_{S_{k-1}} \geq\delta b_n} \pp{(\frac{j}{n},
\frac{X_j}{b_n} )} 
\wc\ppmain
\end{equation}
in $\Mp([0,\infty)\times(0,\infty])$, as $n\rightarrow\infty$, where
$\ppmain$ is defined in \eqref{weenie1}.

\item[(b)]
Suppose additionally that $\bX$ satisfies \eqref{eqecjrvfull} as well
as Conditions \eqref{condeclrd}, \eqref{condeclrdmrv} and
\eqref{condcycmrv}. Recall $\eta^*$ from \eqref{weenie2}.
Then, in $\Mp([0,\infty)\times(0,\infty])$, as $n\to\infty$,
%
%
\begin{equation}
\label{eqppconvxfull} \Nn
\wc
\eta^*,
\end{equation}
and furthermore, the distribution of the cycle maximum has a regularly
varying tail,
%
%
\begin{equation}
\label{eqrvcycfull} t\P_H \Bigl[{b(t)}^{-1}
\sup_{0 \leq j < \ta} X_j \in\cdot \Bigr] \,\mathop{
\longrightarrow}^v\,c \cdot\nua(\cdot) \qquad\mbox{in } \mplus(0,
\infty],
\end{equation}
where
%
%
\begin{equation}
\label{eqrvcycmaxconst} c = \P \Bigl[ \sup_{j \geq1} \xi(j) \leq1 \Bigr] +
\EP \Bigl[\sup_{j \geq1} \xi(j)^\alpha\ind{\sup_{j \geq1} \xi(j)
> 1} \Bigr].
\end{equation}
\end{enumerate}
\end{them}

\begin{pf}
(a) First, note that $\Nnd= \chi_n^0+ \ppmcd$, where
\[
\ppmcd= \sum_{k=0}^\infty\sum
_{j=0}^{\tau^A_{k+1}} \pp{(\frac
{S_k+j}{n},
\frac{X_{S_k+j}}{b_n})} \ind{X_{S_{k}} \geq\delta b_n}.
\]
Hence, by Lemma~\ref{propppneg} it remains to show that $\ppmcd\wc
\ppmain$. Split $\ppmcd$ according to the times $\{\tau_k(b_n)\}$:
\begin{eqnarray*}
\ppmcd&=& \sum_{k=0}^\infty\sum
_{j=0}^{\tkob-1} \pp{(\frac
{S_k+j}{n},
\frac{X_{S_k + j}}{b_n} )} \ind{X_{S_k} \geq\delta b_n} +
\sum_{k=0}^\infty\sum
_{j=\tkob}^{\tau^A_{k+1}} \pp{(\frac
{S_k+j}{n},
\frac{X_{S_k + j}}{b_n} )} \ind{X_{S_k} \geq\delta b_n}\\
& = &\chi'_n+ \chi''_n.
\end{eqnarray*}
The equality holds on the set $\{\tkb\leq\tau^A_k; n\geq1, k\geq1\}$,
which has probability 1 by \eqref{eqcycletimes}.
Because of \eqref{eqcycdistrib} and the independence of the
$(C_k(t),\tkt)$, we have $\chi'_n\,\mathop{=}\limits^{d}\,\ppnmain
$ for each $n$, and
$\ppnmain\wc\ppmain$ by Proposition~\ref{thmppconv}(a).
By Slutsky's theorem, the result follows if
$ \chi''_n\Rightarrow0$, so we show
\[
\P \bigl[ \chi''_n(f) > \gamma \bigr] =
\P \Biggl[ \sum_{k=0}^\infty\sum
_{j=\tkob}^{\tau^A_{k+1}} f \biggl(\frac{S_k+j}{n},
\frac{X_{S_k + j}}{b_n} \biggr) \ind{X_{S_k} \geq\delta b_n} >
\gamma \Biggr] \conv0,
\]
for any $f\in\mathcal{C}^{+}_K([0,\infty)\times(0,\infty])$ and
$\gamma
>0$. Let $f$ have support in $[0,R]\times[M,\infty]$ for integers
$R,M$. The previous
probability is bounded by
\begin{eqnarray*}
&&\P \Biggl[ \sum_{k=0}^{2Rn-1} \sum
_{j=\tkob}^{\tau^A_{k+1}} f \biggl(\frac{S_k+j}{n},
\frac{X_{S_k + j}}{b_n} \biggr) \ind{X_{S_k} \geq\delta b_n} > 0
\Biggr]
\\
&&\quad{} + \P \Biggl[ \sum_{k=2Rn}^\infty\sum
_{j=\tkob}^{\tau^A_{k+1}} f \biggl(\frac{S_k+j}{n},
\frac{X_{S_k + j}}{b_n} \biggr) \ind{X_{S_k} \geq\delta b_n} > 0
\Biggr].
\end{eqnarray*}
Observe that the second term is at most $\P[S_{2Rn}/n \leq R ] = \P
[{S_{2Rn}}/{2Rn} \leq1/2 ] \rightarrow0$ as $n\rightarrow\infty$, since
$S_n/n \rightarrow q$ a.s., and $q\geq1$ by \eqref{eqdefq}.
The first term is bounded by
\begin{eqnarray*}
&&\P \Biggl[ \bigcup_{k=0}^{2Rn-1} \Biggl(
\biggl\{\frac{X_{S_k}}{b_n} \geq\delta \biggr\} \cap\bigcup
_{j=\tkob}^{\tau^A_{k+1}} \biggl\{ \frac
{X_{S_k + j}}{b_n}\geq M \biggr
\} \Biggr) \Biggr]\\
&&\quad{} \leq2R n \P_H \biggl[ \frac{X_0}{b_n} \geq
\delta, \sup_{\tbn< j <
\ta} \frac{X_j}{b_n} \geq M \biggr],
\end{eqnarray*}
which vanishes as $n\rightarrow\infty$ by Condition \eqref{condcyc}.

(b) Recalling the decomposition \eqref{eqnndecomp}, by Lemma~\ref
{propppneg} it is sufficient to show that $\chi^*_n\wc\eta^*$. This
follows by a similar argument as in part (a). Write
\[
\chi^*_n = \sum_{k=0}^\infty\sum
_{j=0}^{\tkob-1} \pp{(\frac
{S_k+j}{n},
\frac{X_{S_k + j}}{b_n} )} + \sum_{k=0}^\infty
\sum_{j=\tkob}^{\tau^A_{k+1}} \pp{(
\frac
{S_k+j}{n}, \frac{X_{S_k + j}}{b_n} )} = {\chi^*_n}'+
{\chi^*_n}''.
\]
Then ${\chi^*_n}'\,\mathop{=}\limits^{d}\,\eta^*_n\wc\eta^*$ by
Proposition~\ref{thmppconv}(b), and Condition \eqref{condcycmrv}
implies that
${\chi^*_n}''\Rightarrow0$.

Next, we show \eqref{eqrvcycfull}.
In light of \eqref{eqcycletimes}, we have
\[
0 \leq t\P_H \biggl[\sup_{0 \leq j < \ta} \frac{X_j}{b(t)} > x
\biggr] - t\P_H \biggl[\sup_{0 \leq j < \tau(b(t))} \frac{X_j}{b(t)} > x
\biggr] \leq t\P_H \biggl[\sup_{\tau(b(t))< j < \ta} \frac{X_j}{b(t)} > x
\biggr] \conv0
\]
under Condition \eqref{condcycmrv}. Recalling that
\[
t\P_H \biggl[\sup_{0 \leq j < \tau(b(t))} \frac{X_j}{b(t)} > x \biggr]
\conv c {x}^{-\alpha}
\]
as $t\rightarrow\infty$ by Proposition~\ref{proprvecmax} (p.
\pageref
{proprvecmax}), where $c$ is as in \eqref{eqrvcycmaxconst}, completes
the proof.
\end{pf}

Setting $M_n = \bigvee_{0 \leq j \leq n} X_j$, Rootz\'en shows
\cite[Theorem~3.2]{rootzen1988maxima} that \eqref{eqrvcycfull}
implies
\[
\P[M_n \leq b_n x] \conv\exp \bigl(-cq^{-1}{x}^{-\alpha}
\bigr),\qquad x>0,
\]
where $c$ is given by \eqref{eqrvcycmaxconst}, and $q$ is the mean
interarrival time \eqref{eqdefq}.
Hence, in the stationary case, $\theta= c/q$ is the extremal index of
the process $\bX$ (\cite[Section~2.2]{leadbetter1988extremal}, \cite
{leadbetter1983extremes}).
On the other hand, for stationary regularly varying Markov chains with
$K \in D(G)$ satisfying a condition analogous to Condition
\eqref{condeclrd}, it is known \cite[Remark~4.7]{basrak2009regularly},
\[
\theta= \P \Bigl[\sup_{j \geq1} Y\xi(j) \leq1 \Bigr] = \P \Bigl[
\sup_{j \geq1} \xi(j) \leq1 \Bigr] - \EP \Bigl[\sup_{j
\geq1}
\xi(j)^\alpha\ind{\sup_{j \geq1} \xi(j) \leq1} \Bigr] 
= c - \EP \Bigl(\sup_{j \geq1} \xi(j)^\alpha \Bigr),
\]
where $Y\sim\mathrm{Pareto}(\alpha)$ supported on $[1,\infty)$,
independent of $\{\xi(j)\}$.
Hence, for a stationary Markov chain $\bX$ satisfying the assumptions
of Theorem~\ref{thmppconvmc}(b), the extremal index is given by
\[
\theta= \frac{1}{q-1} \EP \Bigl(\sup_{j \geq1} \xi(j)^\alpha
\Bigr) = 
\frac{\EP(\sup_{j \geq1} \xi(j)^\alpha)}{\EP_H \ta}.
\]

\subsection{Discussion of conditions}\label{subsec:discussBaby}

We now consider simplifications of the above conditions.

\subsubsection{Cases where $G(\{0\})=0$}\label{subsubsec:case0}

If $G(\{0\})=0$,
we can replace $\bX^{(b(t))}$ with $\bX$ in the finite-dimensional
convergence \eqref{eqecjrv} when $H$ has a regularly varying
tail, meaning that the tail chain approximation completely describes
the extremes of the chain $\bX$ in a finite dimensional sense.
However, $G(\{0\})=0$ also implies that for any $m>0$, as
$t\rightarrow\infty$,
%
%
\begin{equation}
\label{eqcyctimesgz} \P_t\bigl[m < \tat\leq\ta\bigr] \conv1
\end{equation}
(see \cite[Proposition~5.1(d)]{resnick2011asymptotics}) meaning that,
as the initial observation
becomes more extreme, it takes longer for $\bX$ to return to $A$ to
complete the cycle.
Hence, for Condition \eqref{condeclrd} to hold, we need a condition
that ensures
$\bX$ eventually drifts away from extreme states:
%
%
\begin{equation}
\label{condsufflrdz} \lim_{m \rightarrow\infty} \limsup_{t\rightarrow\infty}
\P_{t} \Bigl[ \sup_{m \leq j < \ta} {X_j} > ta \Bigr] = 0
\qquad\mbox{for all } a>0.
\end{equation}

\begin{prop} \label{propsufflrdz}
Suppose $\bX\sim K \in D(G)$ with $G(\{0\})=0$ and positive recurrent
bounded atom $A$, and $X_0 \sim H$ with $1-H \in\RV_{-\alpha}$.
If $\bX$ satisfies Condition \eqref{condsufflrdz}, both Conditions
\eqref{condeclrd} and \eqref{condcyc} hold and
consequently, the convergence \eqref{eqppconvxdfull} takes place.
\end{prop}

\begin{pf}
We first show that
%
%
\begin{equation}
\label{eqpflrdz} \lim_{m \rightarrow\infty} \limsup_{t\rightarrow
\infty} t\P_H
\Bigl[ {X_0}/{b(t)} > \delta, \sup_{m \leq j < {\tau_A}}
{X_j}/{b(t)} > a \Bigr] = 0 \qquad\mbox{for all } a, \delta>0.
\end{equation}
Indeed, for $c > \delta$, we have,
\begin{eqnarray*}
t\P_H \biggl[ \frac{X_0}{b(t)} > \delta, \sup_{m \leq j < {\tau_A}}
\frac{X_j}{b(t)} > a \biggr] &\leq&\int_{[\delta, c]} t
\P_H \biggl[ \frac{X_0}{b(t)} \in du \biggr] \P_{b(t)u}
\biggl[ \sup_{m \leq j < {\tau_A}} \frac{X_j}{b(t)} > a \biggr]
\\
&&{}+ t\P_H \biggl[\frac{X_0}{ b(t)}> c \biggr].
\end{eqnarray*}
Furthermore, for $\delta\leq u \leq c$,
\begin{eqnarray*}
\P_{b(t)u} \biggl[ \sup_{m \leq j < {\tau_A}} \frac{X_j}{b(t)} > a \biggr]
\leq\P_{b(t)u} \biggl[ \sup_{m \leq j < {\tau_A}} \frac
{X_j}{b(t)u} >
\frac{a}{c} \biggr] \leq\sup_{s \geq b(t)\delta} \P_{s} \biggl[
\sup_{m \leq j < {\tau
_A}} \frac{X_j}{s} > \frac{a}{c} \biggr].
\end{eqnarray*}
Hence, by Condition \eqref{condsufflrdz},
\[
\lim_{m \rightarrow\infty} \limsup_{t\rightarrow\infty} t\P_H \biggl[
\frac{X_0}{b(t)} > \delta, \sup_{m \leq j < {\tau_A}} \frac
{X_j}{b(t)} > a \biggr]
\leq\nua[\delta,c] \cdot0 + \nua(c,\infty] = {c}^{-\alpha}.
\]
Letting $c\rightarrow\infty$ establishes \eqref{eqpflrdz}.
As \eqref{eqcycletimes} implies that $\sup_{m \leq j < \tau(b(t))}
X_j \leq\sup_{m \leq j < {\tau_A}} X_j$, Condition \eqref{condeclrd}
follows.
To verify Condition \eqref{condcyc}, argue that
\begin{eqnarray*}
\P_H \Biggl[\frac{X_0}{b(t)} > \delta,\bigvee
_{j=\tau(b(t)) }^ {\tau_A
-1} \frac{{X_j}}{b(t)} > a \Biggr] &\leq&
\P_H \biggl[\frac{X_0}{b(t)} > \delta,\tau \bigl(b(t) \bigr) \leq
m-1 \biggr]
\\
&&{}+ \P_H \Biggl[ \frac{X_0}{b(t)} > \delta, \bigvee
_{j=m}^{ \tau_A-1} \frac{X_j}{b(t)} > a \Biggr],
\end{eqnarray*}
of which the first term vanishes as $t\rightarrow\infty$ because of
\eqref{eqcyctimesgz} (see \cite
[Proposition~5.1(d)]{resnick2011asymptotics}). Appeal to \eqref{eqpflrdz}
and let $m\rightarrow\infty$ to complete the proof.
\end{pf}

Condition \eqref{condsufflrdz} is a condition on the transition kernel
$K$; this is best discussed
by recalling (see \cite[p. 5]{resnick2011asymptotics} for discussion)
that a transition kernel $K\in D(G)$ has an update
function $\psi$ of the form
%
%
\begin{equation}
\label{equpdfunc} \psi \bigl(x, (Z,W) \bigr) = Zx + \phi(x,W),
\end{equation}
where $Z \sim G$ and $t^{-1} \phi(t,w) \rightarrow0$ for $w \in C$ with
$\P[W \in C]=1$ and we
can represent $K$ as
\[
K(x,B) = \P \bigl[\psi \bigl(x,(Z,W) \bigr) \in B \bigr].
\]
Take $V_r = (Z_r,W_r)$, i.i.d. copies of $V=(Z,W)$, and write $\bV_r =
(V_1,\ldots,V_r)$. For $r\geq1$ let $\psi^r(x, \rvect{V}_r)$ denote
the $r$-step update function, i.e., $K^{r}(x,B) = \P[\psi^r(x,\bV_r) \in B]$, and $\psi^0(x) = x$.
By iteration,
\[
\psi^r(x, \bV_r) = \Biggl(\prod
_{j=1}^r Z_j \Biggr) x + \sum
_{\ell= 1}^{r-1} \Biggl( \prod
_{j=\ell+ 1}^r Z_j \Biggr) \phi \bigl(
\psi^{\ell-1}(x,\bV_{\ell-1}),W_\ell \bigr) + \phi \bigl(
\psi^{r-1}(x,\bV_{r-1}), W_r \bigr).
\]
Thus Condition \eqref{condsufflrdz} requires both $Z_m \rightarrow
0$ as in \eqref{eqtchtransient}, and also
an asymptotic stochastic boundedness condition on
$\phi(\cdot, W)$.
Alternately, one could give criteria for Condition \eqref{condsufflrdz}
using mean drift conditions for $\bX$ or $\log\bX$
\cite[p. 229]{meyn:tweedie:1993}.

\subsubsection{Cases where $G(\{0\})> 0$} \label{subsubsec:casePos}

In this case, \eqref{eqecfddconv} implies
\[
\P_{tu} \bigl[\tau(t) = m \bigr] \conv\P \bigl[\tau^* = m \bigr],
\qquad m\geq1,
\]
where $\tau^*$ is a Geometric random variable with parameter $G(\{0\})$.
Hence, the tail chain terminates after a finite number of steps.
If either $y_0(t) \equiv0$ is an extremal boundary, or $K$ satisfies
the regularity condition \eqref{eqregcond} (p. \pageref{eqregcond}),
Theorem~\ref{thmecjrv} assures us that
convergence
\eqref{eqecjrv}
holds for $\bX$ with respect to $\P_H$, and Condition
\eqref{condeclrd} follows directly since
\begin{eqnarray*}
\limsup_{t\to\infty}t\P \Bigl[\sup_{j\geq m} X^{(b(t)}_j/b(t)>a,
X_0>\delta b(t) \Bigr] &\leq&\limsup_{t\to\infty}t \P \bigl[
X_0>\delta b(t), \tau \bigl(b(t) \bigr)\geq m \bigr]
\\
&=& \int_\delta^\infty\nu_\alpha(dx)\P
\bigl[x \xi(m)>0 \bigr] \to0\qquad(m\to\infty).
\end{eqnarray*}

The regularity condition \eqref{eqregcond} extends to any finite
number of steps; that is, iterates of $K$ also satisfy the condition.
However, unless $y_0(t) \equiv0$ is an extremal boundary, we need the
regularity condition to hold uniformly over the whole cycle of random
length $\ta$ to prevent $\bX$ from returning to an extreme state
within the same cycle, after crossing below the extremal boundary.
Condition \eqref{condcycreg} given next accomplishes this.
(Note that even if $y_0(t)\equiv0$ is an extremal boundary for $K$,
we are using an extremal boundary $y(t)$ chosen to satisfy
\eqref{eqcycletimes}.)

%
%
\begin{equation}
\label{condcycreg} \lim_{t\rightarrow\infty} \P_{tu_t} \Bigl[
\sup_{1 \leq j < \ta} {X_j} > {t}a \Bigr] = 0 \qquad\mbox{whenever }
u_t = u(t) \rightarrow0,\qquad a>0.
\end{equation}

Recalling the update function form \eqref{equpdfunc},
the regularity
condition \eqref{eqregcond}
holds if the function $\phi(\cdot,w)$ is bounded near 0 for each $w$
in a set of probability 1
\cite[Proposition~4.1]{resnick2011asymptotics}.
Condition
\eqref{condcycreg} is a stronger boundedness restriction on
$\phi(\cdot,w)$ near 0.
Alternatively, when $K$ satisfies the regularity condition
\eqref{eqregcond}, Condition \eqref{condcycreg} may be viewed as a
restriction on $\ta$, since then
%
%
\begin{equation}
\label{psychoT} \lim_{m\rightarrow\infty} \limsup_{t\rightarrow
\infty} \P_{tu_t}[
\ta> m] = 0 \qquad\mbox{whenever } u_t = u(t) \rightarrow0,
\end{equation}
it is sufficient for \eqref{condcycreg}. This follows from the
decomposition
\[
\P_{tu_t} \Bigl[ \sup_{1 \leq j < \ta} {X_j} > {t}a \Bigr]
\leq\P_{tu_t} [ \tau_A > m ] + \P_{tu_t} \Bigl[
\sup_{1 \leq j \leq m} {X_j} > {t}a \Bigr],
\]
with \eqref{psychoT} controlling the first right-hand term and
\eqref{eqregcond}
controlling the second.

\begin{prop}
Suppose
$\bX\sim K \in D(G)$
with $G(\{0\})>0$
and $1-H \in\RV_{-\alpha}$ and $\bX$
has a positive recurrent, bounded atom $A$.
Then $\bX$ satisfies 
\eqref{condeclrd} with respect to $\P_H$.
Moreover, if either
\begin{enumerate}
\item[(i)] $y_0(t) \equiv0$ is an extremal boundary for $K$,
\end{enumerate}
or
\begin{enumerate}
\item[(ii)] Condition \eqref{condcycreg} holds,
\end{enumerate}
then Condition
\eqref{condcyc} holds with respect to $\P_H$ and thus
convergence \eqref{eqppconvxdfull} takes place.
\end{prop}

\begin{pf}
First, note that by \cite[Proposition~5.1(d)]{resnick2011asymptotics},
as $t\rightarrow\infty$,
%
%
\begin{eqnarray}
\label{eqpfgzpos} t\P_H \biggl[\frac{X_0}{b(t)} > \delta,
\sup_{j \geq m} \frac{{X_j}^{(b(t))}}{b(t)} > a \biggr] &\leq& t\P_H
\biggl[ \frac{X_0}{b(t)} > \delta, \tau\bigl(b(t)\bigr)> m \biggr]
\nonumber
\\[-8pt]
\\[-8pt]
&\to&{\delta}^{-\alpha} \bigl(1-G\bigl(\{0\}\bigr) \bigr)^m.
\nonumber
\end{eqnarray}
Since $G(\{0\})>0$, the right side of \eqref{eqpfgzpos} vanishes as
$m\rightarrow\infty$,
establishing Condition \eqref{condeclrd}.
Next, to analyze Condition
\eqref{condcyc},
consider the case where $y_0(t) \equiv0$ is an extremal
boundary, and write $\tau_0 = \inf\{n \geq0 : X_n = 0\}$. For any
$m$,
%
%
\begin{eqnarray}
\label{eqn:merd} t\P_H \Bigl[{X_0}/{b(t)} > \delta,
\sup_{\tau(b(t))< j < \ta} {{X_j}}/{b(t)} > a \Bigr] &\leq&\sum
_{r=1}^m t \P_H \bigl[{X_0}/{b(t)}
> \delta, \tau\bigl(b(t)\bigr)= r, \tau_0 > r \bigr]
\nonumber
\\
&&{}+ t\P_H \bigl[{X_0}/{b(t)} > \delta, \tau\bigl(b(t)\bigr)> m
\bigr],
\end{eqnarray}
which is obtained by splitting according to whether $\tau(b(t))\leq m$
or the
complement and using the fact that $\tau(b(t))=r$ and
$\sup_{\tau(b(t))< j < \ta} {{X_j}}/{b(t)} > a$ implies
$\tau_0 >r$.
For a typical term in the sum,
\begin{eqnarray*}
t\P_H \bigl[{X_0}/{b(t)} > \delta, \tau\bigl(b(t)\bigr)= r,
\tau_0 > r \bigr] 
&\leq& t\P_H
\bigl[{X_0}/{b(t)} > \delta,\tau_0 > r \bigr]
\\
&&{}- t \P_H \bigl[{X_0}/{b(t)} > \delta, \tau\bigl(b(t)\bigr)> r
\bigr]
\\
&\conv&{\delta}^{-\alpha} \bigl(1-G\bigl(\{0\}\bigr)\bigr)^r -
{\delta}^{-\alpha
} \bigl(1-G\bigl(\{0\}\bigr)\bigr)^r = 0,
\end{eqnarray*}
the convergence following from \cite
[Proposition~5.1(d)]{resnick2011asymptotics},
since both $y(t)$ and $y_0 (t)$ are extremal boundaries.
The right most term in \eqref{eqn:merd} is handled as in \eqref{eqpfgzpos}.

Finally analyze Condition \eqref{condcyc} when
Condition \eqref{condcycreg} holds.
For any $m$, we have
\begin{eqnarray*}
&&t\P_H \Bigl[{X_0}/{b(t)} > \delta,
\sup_{\tau(b(t))< j < \ta} {{X_j}}/{b(t)} > a \Bigr]
\\
&&\quad\leq\sum_{r=1}^m t
\P_H \Bigl[{X_0}/{b(t)} > \delta, \sup_{r < j <
\ta}
{{X_j}}/{b(t)} > a, \tau\bigl(b(t)\bigr)= r \Bigr]
\\
&&\phantom{\quad\leq} {}+ t\P_H \bigl[{X_0}/{b(t)} >
\delta, \tau\bigl(b(t)\bigr)> m \bigr],
\end{eqnarray*}
and
%
%
\begin{eqnarray}
\label{eqpfintht} &&t\P_H \Bigl[{X_0}/{b(t)} > \delta,
\sup_{r < j < \ta} {{X_j}}/{b(t)} > a, \tau\bigl(b(t)\bigr)= r \Bigr]
\\
&&\quad= \int_{(\delta,\infty] \times(y(b(t)),\infty]^{r-1}
\times
[0,\infty]} t\P_H \bigl[ (
{X_0}, {\bX_r} )/b(t) \in d(x_0,\boldsymbol
x_r) \bigr] h_t(x_r),
\nonumber
\end{eqnarray}
where
\[
h_t(x) = \ind{ [0,y (b(t) ) ]} (x)
\P_{b(t)x} \Bigl[ \sup_{1 \leq j <
\ta} {X_j}/{b(t)} > a
\Bigr].
\]
We claim that $h_t(u_t) \rightarrow0$ whenever $u_t \rightarrow u \geq
0$. Indeed, if $u>0$, then $h_t(u_t) = 0$ for large $t$ such that
$y(b(t))<u$.
Otherwise, $u_t\rightarrow0$, and $h_t(u_t) \rightarrow0$ by
Condition \eqref{condcycreg}. Therefore, the integral converges to 0 by
combining Lemmas~8.2 and~8.4 with Theorem~3.2 from~\cite{resnick2011asymptotics}.
Applying \eqref{eqpfgzpos} completes the proof.
\end{pf}

\subsection{Weak convergence to a cluster process}\label{subsec:weak}

If the finite-dimensional distributions of ${\bX}$ are jointly
regularly varying (in the sense of \eqref{eqecjrvfull} with $\bX$
replacing $\bX^{(b(t))}$), we obtain a point
process limit for $\bX$ under a condition
analogous to Condition \eqref{condeclrdmrv}: \textit{There exists
$m'_0\geq1$ such that}
%
%
\begin{equation}
\label{condsufflrdfull} \lim_{\delta\downarrow0} \limsup_{t\rightarrow\infty} t
\P_H \Bigl[{X_0}/{b(t)} \leq\delta, \sup_{m_0' \leq j < \ta}
{X_j}/{b(t)} > a \Bigr] = 0 \qquad\mbox{for all } a>0.
\end{equation}

\begin{prop}
Suppose $\bX\sim K \in D(G)$ has a positive recurrent, bounded atom
$A$, and
$1-H \in\RV_{-\alpha}$. Assume further that, with respect to $\P_H$, $\bX$ is regularly varying in the sense of
\eqref{eqecjrvfull}, with $\bX$ replacing $\bX^{(b(t))}$, and
satisfies Condition \eqref{condcyc}.
Under Condition \eqref{condsufflrdfull}, both Conditions \eqref
{condeclrdmrv} and \eqref{condcycmrv} hold with respect to $\P_H$.
\end{prop}

\begin{pf}
Recalling $\sup_{m \leq j < \tau(b(t))} X_j \leq\sup_{m \leq j < \ta} X_j$
yields \eqref{condeclrdmrv}. Next, given $\delta>0$, write
\begin{eqnarray*}
\P_H \biggl[ \sup_{\tau(b(t))< j < \ta} \frac{{X_j}}{b(t)} > a \biggr] &
\leq&\P_H \biggl[\frac{X_0}{b(t)} > \delta, \sup_{\tau(b(t))< j < \ta}
\frac{{X_j}}{b(t)} > a \biggr]
\\
&&{}+ \P_H \biggl[\frac{X_0}{b(t)} \leq \delta,
\sup_{1 \leq j < \ta} \frac{{X_j}}{b(t)} > a \biggr].
\end{eqnarray*}
Condition \eqref{condcyc} makes the first right side term go to $0$
as $t \to\infty$ and for the second term we have,
\begin{eqnarray*}
&&\limsup_{t\rightarrow\infty} \P_H \biggl[\frac{X_0}{b(t)} \leq \delta
, \sup_{1 \leq j < \ta} \frac{{X_j}}{b(t)} > a \biggr]
\\
&&\quad\leq\limsup_{t\rightarrow\infty} \P_H \biggl[\frac
{X_0}{b(t)}
\leq \delta, \sup_{1 \leq j < m_0'} \frac{{X_j}}{b(t)} > a \biggr] 
+
\limsup_{t\rightarrow\infty} \P_H \biggl[\frac{X_0}{b(t)} \leq \delta,
\sup_{m_0' \leq j < \ta} \frac{{X_j}}{b(t)} > a \biggr]
\\
&&\quad= \mu^* \bigl([0,\delta] \times{[\bz,\boldsymbol {a}]}^{\mathrm{c}} \bigr) +
\limsup_{t\rightarrow\infty} \P_H \biggl[\frac{X_0}{b(t)} \leq \delta,
\sup_{m_0' \leq j < \ta} \frac{{X_j}}{b(t)} > a \biggr]
\end{eqnarray*}
where $\boldsymbol{a} = (a,\ldots,a)$. Letting $\delta\downarrow0$,
the first term vanishes by \eqref{eqecjrvfull},
since $\mu^*(\E_\sqsupset^*\setminus\E_\sqsupset)=0$ \cite
[Theorem~5.1]{resnick2011asymptotics}. The second term is taken care of by
Condition \eqref{condsufflrdfull}.
\end{pf}

We now rephrase Theorem~\ref{thmppconvmc} in terms of our new conditions.

\begin{them} \label{thmppconvg}
Let $\bX$ be a Markov chain on $[0,\infty)$ with transition kernel
$K\in D(G)$ such that $K$ has a positive recurrent bounded atom in the
sense of \eqref{eqatom}, \eqref{eqposrec}, and \eqref{eqatombdd}. The
initial distribution $H$ has a regularly varying tail and
satisfies $tH(b(t)\cdot) \,\mathop{\rightarrow}\limits^v\,\nu_\alpha(\cdot)$ in $\mplus
(0,\infty]$.
Assume for any $m\geq0$ that $(X_0,\ldots,X_m)$ is regularly varying
in $\mplus([0,\infty]^{m+1} \setminus\{\bz\})$,
%
%
\begin{equation}
\label{eqjrvfull} t\P_H \bigl[ ({{X_0}}, \ldots,
{{X_m}} )/b(t) \in(dx_0, d\boldsymbol x_m)
\bigr] \,\mathop{\longrightarrow}^v\,\nu_\alpha(dx_0)
\P_{x_0} \bigl[ (T_1,\ldots, T_m) \in d
\boldsymbol x_m \bigr]
\end{equation}
and that Condition \eqref{condsufflrdfull} holds with respect to $\P_H$.
\begin{enumerate}[(b)]
\item[(a)] If $G(\{0\})= 0$, and $K$ satisfies Condition \eqref
{condsufflrdz}, then
\[
\sum_{j=0}^{\infty} \pp{(
\frac{j}{n}, \frac{X_j}{b_n})} 
\wc\sum
_{k=0}^\infty\sum_{j=0}^{\infty}
\pp{(qt_k, i_k\xikj)} \qquad\mbox{in } \Mp \bigl([0,
\infty )\times (0,\infty] \bigr) \qquad\mbox{as } n\rightarrow\infty.
\]

\item[(b)] If $G(\{0\})> 0$, and either $y_0(t) \equiv0$ is an extremal
boundary for $K$, or $K$ satisfies Condition \eqref{condcycreg}, then
\[
\sum_{j=0}^{\infty} \pp{(
\frac{j}{n}, \frac{X_j}{b_n} )} 
\wc\sum
_{k=0}^\infty\sum_{j=0}^{\tau^*_{k+1}-1}
\pp{(qt_k, i_k\xikj)} \qquad\mbox{in } \Mp \bigl([0,
\infty )\times (0,\infty] \bigr) \qquad\mbox{as } n\rightarrow\infty,
\]
where the $\{\tau^*_k\}$ are i.i.d. Geometric random variables with
parameter $G(\{0\})$.
\end{enumerate}
\end{them}

%
%
%
%
%


\section{\texorpdfstring{Proof of Proposition \protect\ref{thmppconv}}{Proof of Proposition 3.1}}
\label{secppconv}


Recall that $\{\bX,\bX_k,k\geq0\}$ are i.i.d. copies of a Markov chain
$\bX\sim K$ with heavy tailed initial distribution $H$ satisfying
$tH((b(t) \cdot)
\,\mathop{\rightarrow}\limits^v\,\nu_\alpha(\cdot)$ in $\mathbb
{M}_+(0,\infty]$
where $b(t)\rightarrow\infty$. The extremal boundary downcrossing
time by
$\bX_k$ is $\{\tkit\}$.
Let
$\{\bxi,\bxi_k,k\geq0\}$ be i.i.d. copies of the multiplicative random
walk $\bxi=\{\xi(m),m\geq0\}$.
The hitting time of $0$ by $\bxi_k$ is $\tau^*_k$.
A~PRM on $[0,\infty)\times(0,\infty]$ with
mean measure $\leb\times\nua$, independent of the
$\{\bxi,\bxi_k,k\geq0\}$ is
$\limpp=
\sum\pp{(t_k, i_k)}$ and
$\{S_k\}$ is the renewal process given by \eqref{eqrenewaltimes} with
finite mean interarrival time $q$.
For convenience, write $\bX_{k,m}^{(t)} = (\Xko^{(t)},\ldots, \Xkm^{(t)})$ and $\xivkm= (\xiko,\ldots,\xikm)$.

\begin{pf*}{Proof of Proposition~\ref{thmppconv}}
(a) First, recall that, under our assumptions $K\in D(G)$ and $H$
having a
regularly varying tail, the convergence
\eqref{eqecjrv} takes place for the chain $\bX$ on the space $\E_\sqsupset:=
(0,\infty] \times[0,\infty]^m$, with limit measure $\mu$ given by
\eqref{eqlimmeas}. This implies that
\cite[Corollary~6.1, p.~183]{resnick2007heavy},
%
%
\begin{equation}
\label{eq:prmLim} \sum_{k=0}^\infty\pp{
( {k}/{n}, {\bX_{k,m}^{(b_n)}}/{b_n} )} =
\sum_{k=0}^\infty \pp{( {k}/{n}, (
{\Xko^{(b_n)}},\ldots, {X_{km}^{(b_n)}}
)/{b_n} )} \wc\prm(\leb\times\mu),
\end{equation}
in $\Mp([0,\infty) \times\E_\sqsupset)$.
Since $\sum_{k=0}^\infty\pp{(t_k,i_k,\xivkm)}$
is PRM on $[0,\infty) \times\E_\sqsupset$
(\cite[Proposition~5.3, p.~123]{resnick2007heavy}), a mapping argument
(\cite[Proposition~5.2, p. 121]{resnick2007heavy}) implies
\[
\sum_{k=0}^\infty\pp{(t_k,
i_k\xivkm)} = \sum_{k=0}^\infty
\pp{(t_k, i_k, i_k
\xi_k(1), \ldots, i_k\xikm )} \sim\prm(\leb\times
\mu),
\]
in $ \Mp([0,\infty) \times\E_\sqsupset)$,
by \eqref{eqlimmeas}.
So we can rewrite \eqref{eq:prmLim} in $\Mp( [0,\infty) \times
\E_\sqsupset)$,
%
%
\begin{equation}
\label{eqpfppvecconv} \ppnveco= \sum_{k=0}^\infty
\pp{( {k}/{n}, {\bX_{k,m}^{(b_n)}}/{b_n} )}
\wc\sum_{k=0}^\infty \pp{(t_k,
i_k\xivkm)} = \ppveco.
\end{equation}

Second, we rescale the time axis to place points at the epochs
$S_k$. (See \cite{mikosch2006activity}.)
The counting function for the points $\{S_k\}$ is
$N(t) = \sum_k \pp{S_k} [0,t]$
and $ N^{\leftarrow}(t) = \inf\{s: N(s)\geq t \} = S_{[t]} $ is
the left continuous inverse process.
Define $\Theta_n(\cdot) = n^{-1}N^{\leftarrow}(n \cdot)$, so that
$S_k/n = \Theta_n(k/n)$
and $\Theta_n$ is a random element of $D^{\uparrow}[0,\infty)$, the subspace
of non-decreasing elements of $D_{\mathrm{left}}[0,\infty)$.
By the Strong Law of Large Numbers, with probability $1$,
\[
\Theta_n(t) = \frac{[nt]}{n} \frac{S_{[nt]}}{[nt]} \conv t\cdot q,
\qquad t\geq0,
\]
so $\Theta_n (\cdot) \rightarrow q(\cdot)$ in $D^{\uparrow
}[0,\infty)$.
We transform time points using the mapping $T_1 : D^{\uparrow
}[0,\infty)
\times
\mplus([0,\infty) \times\E_\sqsupset) \mapsto\mplus([0,\infty)
\times\E_\sqsupset)$
given by
%
%
\begin{equation}
\label{eqpftimech} T_1{m}(f) = \iint f \bigl(x(u),v \bigr)
m(du,dv), \qquad f\in\mathcal{C}^{+}_K \bigl([0,\infty
) \times \E_\sqsupset\bigr).
\end{equation}
Applying \cite[Proposition~3.1, p. 57]{resnick2007heavy} to
\eqref{eqpfppvecconv}, we have $(\Theta_n (\cdot),\ppnveco) \wc
(q(\cdot),\ppveco)$ in
$D^{\uparrow}[0,\infty) \times\Mp([0,\infty) \times\E_\sqsupset)$.
Since $T_1$ is a.s. continuous at $(q(\cdot),\ppveco)$ (Lemma~\ref
{lemcontvagueconv}, p. \pageref{lemcontvagueconv}), the Continuous
Mapping Theorem gives in $\Mp([0,\infty) \times\E_\sqsupset)$,
%
%
\begin{equation}
\label{eq:toRestrict} \eta'_n= \sum
_{k=0}^\infty\pp{( {S_k}/{n}, {
\bX_{k,m}^{(b_n)}}/{b_n} )} = T_1(
\Theta_n,\ppnveco) \wc T_1 \bigl(q(\cdot),\ppveco \bigr)
= \sum_{k=0}^\infty \pp{(qt_k,
i_k\xivkm)} = \eta'.
\end{equation}

Now stack the components of $\bX_{k,m}^{(b_n)}$ above the time point $S_k/n$.
To make functionals continuous, it is necessary to compactify the state
space by letting $\Ld
:= [\delta,\infty]\times[0,\infty]^m$.
Define the restriction functional $T_2 : \Mp([0,\infty)
\times\E_\sqsupset) \rightarrow\Mp([0,\infty) \times
\Ld)$ by $T_2 m = m(\cdot\cap([0,\infty) \times\Ld))$.  From
\cite[Proposition~3.3]{feigin1996parameter}, $T_2$ is almost
surely continuous at $\eta'$ provided $\P[ \eta'(\bdry
([0,\infty)\times\Ld)) = 0 ] = 1$ and since $\EP(\eta'(\bdry
([0,\infty)\times\Ld)))=0$ due to $\nu_\alpha(\{\delta\})=0$, the
a.s. continuity is verified.
Therefore, in $\Mp([0,\infty)\times\Ld)$, the restricted version of
\eqref{eq:toRestrict} is $\eta''_n:= T_2(\eta'_n) \wc T_2(\eta')
=:\eta''$.
Define the stacking functional $T_3 : \Mp([0,\infty)\times\Ld)
\rightarrow\Mp([0,\infty)\times[0,\infty])$ by
\[
T_3 \biggl(\sum_k
\epsilon_{ (t_k,y_k(0),\ldots,y_k(m) )} \biggr) =\sum_k \sum
_{j=0}^m\epsilon_{(t_k,y_k(j))}
\]
or for $m\in\Mp([0,\infty)\times\Ld)$, $f\in\mathcal
{C}^{+}_K([0,\infty)
\times[0,\infty])$,
$ T_3{m}(f) = \iint\{ \sum_{j=0}^m f(u,v_j) \} m(du,
d\rvect{v})$.
Given such $f$ with support in $[0,R]\times[0,\infty]$, $R$ a
positive integer, the function
$\varphi(u,\rvect{v}) := \sum_{j=0}^m f(u,v_j) \in\mathcal{C}^{+}_K
([0,\infty)\times\Ld)$, since it is clearly non-negative continuous,
and $\varphi=0$ outside of $[0,R]\times\Ld$.
The continuity of $T_3$ is clear: given $m_n \,\mathop{\rightarrow
}\limits^v\,m$ in
$\Mp([0,\infty)\times\Ld)$, we have $T_3{m}_n(f) = m_n(\varphi)
\rightarrow m(\varphi) = T_3{m}(f)$.
Consequently, in $\Mp([0,\infty)\times[0,\infty])$,
%
%
\begin{equation}
\label{eqwcppnvst} \ppnvst= \sum_{k=0}^\infty
\sum_{j=0}^{m} \pp{(
\frac{S_k}{n}, \frac{\Xkj^{(b_n)}}{b_n})} \ind{\frac{\Xko
}{b_n}\geq\delta}
= T_3 \bigl(\eta''_n \bigr)
\wc T_3 \bigl(\eta'' \bigr) = \sum
_{k=0}^\infty\sum
_{j=0}^{m} \pp{(qt_k, i_k
\xikj)} \ind{i_k \geq\delta} = \ppvst.
\end{equation}

Now adjust the sum over $j$ to replace $X_{kj}^{(b_n)}$ with $X_{kj}$.
From \eqref{eqwcppnvst}, we readily get,
%
%
\begin{equation}
\label{eqpfpptrunc} \ppnvst \bigl( \cdot\cap \bigl([0,\infty )\times
(0, \infty] \bigr) \bigr) \wc\ppvst \bigl( \cdot\cap \bigl([0,\infty )
\times (0,\infty] \bigr) \bigr) \qquad\mbox{in } \Mp \bigl([0,\infty )
\times (0,\infty] \bigr)
\end{equation}
%
by noting that any $f\in\mathcal{C}^{+}_K([0,\infty)\times(0,\infty])$
extends to $\bar{f} \in\mathcal{C}^{+}_K([0,\infty)\times[0,\infty
])$ with
$\bar{f}(s,0)=0$ for $s\geq0$.
Moreover, recalling $\{\tkoit\}$ and $\{\tau^*_{k+1}\}$, the first hitting
times of 0 by $\{\bX^{(t)}_k\}$ and $\{\bxi_k\}$ respectively, put
\[
\skt= \tkoit\smin(m+1) \quad\mbox{and}\quad\sigma^*_{k+1}=
\tau^*_{k+1}\smin(m+1),\qquad k\geq0.
\]
%
Using this notation, the convergence \eqref{eqpfpptrunc} becomes, in
$\Mp([0,\infty)\times(0,\infty])$,
%
%
\begin{equation}
\label{eqwcppnst} \ppnst= \sum_{k=0}^\infty
\sum_{j=0}^{\sk-1} \pp{(
\frac
{S_k}{n}, \frac{\Xkj}{b_n} )} \ind{\frac{\Xko}{b_n}\geq \delta}
\wc\sum_{k=0}^\infty\sum
_{j=0}^{\sigma^*_{k+1}-1} \pp{(qt_k, i_k
\xikj)} \ind{i_k \geq\delta} = \ppst.
\end{equation}
Equation \eqref{eqwcppnst} allows spreading the stacks of $\ppnst$ in
time and we verify
%
%
\begin{equation}
\label{eqwcppnstj} \tilde{\eta}^*_n= \sum
_{k=0}^\infty\sum_{j=0}^{\sk-1}
\pp{(\frac
{S_k+j}{n}, \frac{\Xkj}{b_n})} \ind{\frac{\Xko}{b_n}
\geq\delta} \wc\ppst\qquad\mbox{in } \Mp \bigl([0,\infty )\times (0,
\infty] \bigr).
\end{equation}
This follows from Slutsky's theorem if $d_v(\tilde{\eta}^*_n,\ppnst)
\,\mathop{\rightarrow}\limits^{P}\,0$, where $d_v$ is the vague
metric on $\Mp([0,\infty)\times
(0,\infty])$; it suffices to show that $\P[|\tilde{\eta}^*_n(f) -
\ppnst(f)| > \gamma] \rightarrow0$ for any $f\in\mathcal
{C}^{+}_K([0,\infty
)\times
(0,\infty])$ and $\gamma>0$.
For such $f$ with support in $[0,R] \times[M,\infty]$, for $R,M$
positive integers, we have
\begin{eqnarray*}
&&\P \bigl[ \bigl|\tilde{\eta}^*_n(f) - \ppnst(f)\bigr| > \gamma \bigr]
\\
&&\quad{}= \P \Biggl[ \Biggl|\sum_{k=0}^\infty\sum
_{j=0}^{\sk-1} f \biggl(\frac{S_k+j}{n},
\frac{\Xkj}{b_n} \biggr) \ind{\frac{\Xko
}{b_n}\geq\delta}\\
&&\hphantom{\quad{}= \P \Biggl[ \Biggl|} - \sum
_{k=0}^\infty\sum_{j=0}^{\sk-1}
f \biggl(\frac{S_k}{n}, \frac
{\Xkj}{b_n} \biggr) \ind{\frac{\Xko}{b_n}
\geq\delta} \Biggr| > \gamma \Biggr]
\\
&&\quad{}\leq\P \Biggl[ \sum_{k=0}^\infty\sum
_{j=1}^{\sk-1} \biggl| f \biggl( \frac{S_k+j}{n},
\frac{\Xkj}{b_n} \biggr) - f \biggl( \frac{S_k}{n}, \frac{\Xkj}{b_n}
\biggr) \biggr|
\\
&&\phantom{\quad\leq\P \Biggl[ \sum_{k=0}^\infty\sum
_{j=1}^{\sk-1}} {}\times\epsilon_{ (\frac{S_k}{n}, \frac
{\Xkj}{b_n} ) } \bigl([0,R]\times[M,
\infty] \bigr) \ind{ \frac{\Xko}{b_n}\geq\delta} > \gamma \Biggr].
\end{eqnarray*}
Since $f$ is uniformly continuous, given $\rho>0$, there exists $v>0$
such that $|f(x)-f(y)|<\rho$ whenever $\Vert x-y \Vert< v$.
For $n$ so large that $m/n < v$, we have
\begin{eqnarray*}
&&\sum_{k=0}^\infty\sum
_{j=1}^{\sk-1} \biggl| f \biggl( \frac{S_k+j}{n},
\frac{\Xkj}{b_n} \biggr) - f \biggl( \frac{S_k}{n}, \frac{\Xkj}{b_n}
\biggr)\biggr | \epsilon_{ (\frac{S_k}{n}, \frac{\Xkj}{b_n} ) } \bigl([0,R]\times[M,\infty] \bigr) \ind{
\frac{\Xko}{b_n}\geq\delta}
\\
&&\quad< \rho\cdot\ppnst \bigl([0,R]\times[M,\infty] \bigr),
\end{eqnarray*}
implying that
\begin{eqnarray*}
\limsup_{n\rightarrow\infty} \P \bigl[\bigl|\tilde{\eta}^*_n(f) - \ppnst(f)\bigr| >
\gamma \bigr] &\leq&\limsup_{n\rightarrow\infty} \P \bigl[\ppnst \bigl([0,R]\times[M,
\infty] \bigr) \geq\gamma\rho^{-1} \bigr]
\\
&\leq&\P \bigl[\ppst \bigl([0,R]\times[M,\infty] \bigr) \geq \gamma
\rho^{-1} \bigr]
\end{eqnarray*}
by \eqref{eqwcppnst}.
So \eqref{eqwcppnstj} follows by letting $\rho\rightarrow0$.

Finally, we remove the restriction in $m$ on the stacks. Recall the
definitions of $\eta_n$ and $\eta$ from Proposition~\ref{thmppconv}.
To apply a Slutsky argument (e.g., \cite{resnick2007heavy},
Theorem~3.5, p. 56),
we show, for $\gamma>0$,
%
%
\begin{equation}
\label{eqpfconvtog} \lim_{m\rightarrow\infty} \P \bigl[ d_v(\ppst,
\ppmain) > \gamma \bigr] = 0 \quad\mbox{and}\quad\lim_{m\rightarrow\infty}
\limsup_{n \rightarrow\infty} \P \bigl[ d_v \bigl(\tilde{
\eta}^*_n, \ppnmain \bigr) > \gamma \bigr] = 0.
\end{equation}
Let $f\in\mathcal{C}^{+}_K([0,\infty)\times(0,\infty])$ with support
$[0,R]\times[M,\infty]$.
Taking $\delta<a<\infty$, we write
\[
\bigl|\ppst(f) - \ppmain(f)\bigr| = \sum_{k=0}^\infty
\sum_{j=m+1}^{\infty} f\bigl(qt_k,
i_k\xikj\bigr) \cdot( \ind{\delta\leq i_k < a} +
\ind{i_k \geq a} ).
\]
Hence,
\begin{eqnarray*}
\P \bigl[ \bigl|\ppst(f) - \ppmain(f)\bigr| > \gamma \bigr] &\leq&\P \Biggl[ \sum
_{k=0}^\infty\sum_{j=m+1}^{\infty}
f\bigl(qt_k, i_k\xikj\bigr) \ind{\delta\leq i_k < a}
> \gamma/2 \Biggr]
\\
&&{} + \P \Biggl[ \sum_{k=0}^\infty\sum
_{j=m+1}^{\infty} f\bigl(qt_k,
i_k\xikj\bigr) \ind{i_k \geq a} > \gamma/2 \Biggr]
\\
&=&A+B.
\end{eqnarray*}
Writing $\xi^*_k(m) = \sup_{j\geq m+1}\xikj$ for $k \geq0$,
term $A$ is bounded by
\begin{eqnarray*}
&&\P \Biggl[ \sum_{k=0}^\infty \Biggl\{
\pp{(t_k, i_k)} \bigl( [0,R/q] \times[\delta,a
) \bigr) \sum_{j=m+1}^{\infty} \ind{\xikj> \frac
Ma} \Biggr\}> 0 \Biggr]
\\
&&\quad\leq\P \bigl[\zeta'_m \bigl([0,R/q] \times[
\delta,\infty] \times(M/{a}, \infty) \bigr) > 0 \bigr],
\end{eqnarray*}
where, since $\{\bxi_k\}$ are i.i.d. and independent of $\zeta$, in
$\Mp
([0,\infty)\times(0,\infty] \times[0,\infty])$,
\[
\zeta'_m = \sum_{k=0}^\infty
\pp{(t_k, i_k, \xi^*_k(m) )}
\sim \prm \Bigl(\leb\times\nua\times\P \Bigl[\sup_{j\geq m+1} \xi(j) \in\cdot
\Bigr] \Bigr).
\]
Therefore,
$\P[\zeta'_m([0,R/q] \times[\delta,\infty] \times(Ma^{-1},
\infty)) > 0] = 1-\exp\{-\lambda\}$, where
\begin{eqnarray*}
\lambda&=& \leb[0,R/q] \cdot\nua[\delta,\infty] \cdot\P \Bigl[
\sup_{j \geq m+1} \xi(j) > Ma^{-1} \Bigr]
\\
&=& R q^{-1} {
\delta}^{-\alpha} \P \Bigl[ \sup_{j \geq m+1} \xi(j) > Ma^{-1}
\Bigr]
\\
&\conv& 0
\end{eqnarray*}
as $m\rightarrow\infty$ by \eqref{eqcondeclrdtch}, a consequence of
Condition \eqref{condeclrd}.
For term $B$, we have the bound
\[
\P \bigl[ \zeta \bigl([0,R/q]\times[a,\infty] \bigr) > 0 \bigr] = 1-\exp \bigl
\{-\EP\zeta \bigl([0,R/q] \times[a,\infty] \bigr) \bigr\} = 1-\exp \bigl
\{-Rq^{-1}{a}^{-\alpha} \bigr\}.
\]
Letting $a\rightarrow\infty$ establishes the first limit in
\eqref{eqpfconvtog}.

To prove the second limit in \eqref{eqpfconvtog}, observe that
\begin{eqnarray*}
\P \bigl[\bigl|\tilde{\eta}^*_n(f) - \ppnmain(f)\bigr| > \gamma \bigr] &=& \P
\Biggl[\sum_{k=0}^\infty\sum
_{j=m+1}^{\tk-1} f \biggl(\frac
{S_k+j}{n},
\frac{\Xkj}{b_n} \biggr) \ind{\frac{\Xko}{b_n}\geq\delta} > \gamma \Biggr]
\\
&\leq&\P \Biggl[\sum_{k=0}^{2Rn-1} \sum
_{j=m+1}^{\infty} f \biggl(\frac{S_k+j}{n},
\frac{\Xkj^{(b_n)}}{b_n} \biggr) \ind{\frac{\Xko}{b_n}\geq\delta } > 0 \Biggr]
\\
&&{} + \P \Biggl[\sum_{k=2Rn}^\infty\sum
_{j=m+1}^{\infty} f \biggl(\frac
{S_k+j}{n},
\frac{\Xkj^{(b_n)}}{b_n} \biggr) \ind{\frac{\Xko}{b_n}\geq\delta } > 0 \Biggr].
\end{eqnarray*}
The first term is bounded by
\[
\P \Biggl[ \bigcup_{k=0}^{2Rn-1} \Biggl(
\biggl\{\frac{\Xko}{b_n} \geq\delta \biggr\} \cap\bigcup
_{j=m+1}^\infty \biggl\{ \frac{\Xkj^{(b_n)}
}{b_n} \geq M \biggr
\} \Biggr) \Biggr] \leq2Rn \P \biggl[\frac{X_0}{b_n} \geq\delta,
\sup_{j \geq m+1} \frac{{X_j}^{(b_n)}}{b_n} \geq M \biggr],
\]
and
\[
\lim_{m\rightarrow\infty} \limsup_{n \rightarrow\infty} n \P \biggl[\frac
{X_0}{b_n} \geq
\delta, \sup_{j \geq m+1} \frac{{X_j}^{(b_n)}}{b_n} \geq M \biggr] = 0
\]
by Condition \eqref{condeclrd}.
The second term is at most $\P[S_{2Rn}/n \leq R ] = \P[{S_{2Rn}}/{2Rn}
\leq1/2 ] \rightarrow0$ as $n\rightarrow\infty$, since $S_n/n
\rightarrow q$
a.s., and $q\geq1$ by \eqref{eqdefq}. This establishes
\eqref{eqpfconvtog}, completing the proof of part (a).

(b) This amounts to removing the restrictions in $\delta$, under
the additional assumptions
\eqref{eqecjrvfull} and Condition \eqref{condeclrdmrv}.
We proceed via a Slutsky argument showing that for any $\gamma>0$,
%
%
\begin{equation}
\label{eqconvtogrv} \lim_{\delta\rightarrow0} \P \bigl[ d_v \bigl(
\ppmain, \eta^* \bigr) > \gamma \bigr] = 0 
\quad\mbox{and}
\quad \lim_{\delta\rightarrow0} \limsup_{n \rightarrow\infty} \P \bigl[ d_v
\bigl( \ppnmain, \eta^*_n \bigr) > \gamma \bigr] = 0.
\end{equation}
%
Let $f\in\mathcal{C}^{+}_K([0,\infty)\times(0,\infty])$ with support
$[0,R]\times[M,\infty]$, and
note that
\[
\bigl|\ppmain(f)-\eta^*(f)\bigr| = \sum_{k=0}^\infty
\sum_{j=0}^{\tau^*_{k+1}
-1} f\bigl(qt_k,
i_k\xikj\bigr) \ind{i_k < \delta}.
\]
Hence, writing $\xi^*_k = \sup_{j\geq1} \xikj$, and $\zeta' =
\sum_{k=0}^\infty\pp{(t_k, i_k, i_k\xi^*_k)}$, we have
\begin{eqnarray*}
\P \bigl[ \bigl|\ppmain(f) - \eta^*(f)\bigr| > \gamma \bigr] &\leq& \P \Biggl[ \sum
_{k=0}^\infty\pp{(t_k,
i_k)} [0,R/q] \times(0,\delta) \sum_{j=0}^{\infty}
\ind{i_k \xikj\geq M} > 0 \Biggr]
\\
&\leq& \P \bigl[ \zeta' \bigl([0,R/q] \times(0,\delta) \times[M,
\infty] \bigr) > 0 \bigr].
\end{eqnarray*}
The $\{\bxi_k\}$ are i.i.d. and independent of $\limpp$, so $\zeta'
\sim\prm(\mu')$ on $\Mp[0,\infty)\times(0,\infty] \times
[0,\infty]$ with
\[
\mu'(ds,dx,dy) = \leb(dx) \cdot\nua(dx) \cdot\P \Bigl[
\sup_{j\geq1} \xi(j) \in x^{-1} dy \Bigr]
\]
by \cite[Proposition~5.6, p. 144]{resnick2007heavy}.
Therefore, $\P[ \zeta' ([0,R/q] \times(0,\delta) \times
[M,\infty]) > 0 ] = 1-\exp\{-\lambda\}$, where by Lemma~\ref
{lemintmom},
\[
\lambda= R q^{-1} \int_{(0,\delta)} \nua(dx) \P \Bigl[
\sup_{j\geq1} \xi(j) \geq Mx^{-1} \Bigr] \leq R
q^{-1}{M}^{-\alpha} \EP \Bigl[ \sup_{j\geq1}
\xi(j)^\alpha\cdot\ind{ \sup_{j\geq1} \xi(j) > M
\delta^{-1}} \Bigr].
\]
Apply \eqref{eqcondeclrdtchmom} and dominated convergence
as $\delta\downarrow0$ to get
$\lambda\rightarrow0$ and hence the first
limit in \eqref{eqconvtogrv}.

For the second limit in \eqref{eqconvtogrv}, we have
\begin{eqnarray*}
\P \bigl[\bigl|\ppnmain(f) - \eta^*_n(f)\bigr| > \gamma \bigr] &=& \P \Biggl[
\sum_{k=0}^\infty\sum
_{j=0}^{\tk-1} f \biggl(\frac
{S_k+j}{n},
\frac{\Xkj}{b_n} \biggr) \ind{\frac{\Xko^{(b_n)}}{b_n} <\delta} > \gamma \Biggr]
\\
&\leq&\P \Biggl[\sum_{k=0}^{2Rn-1} \sum
_{j=0}^{\tk-1} f \biggl(\frac{S_k+j}{n},
\frac{\Xkj}{b_n} \biggr) \ind{\frac{\Xko^{(b_n)}}{b_n} < \delta} > 0 \Biggr]
\\
&&{} + \P \Biggl[\sum_{k=2Rn}^\infty\sum
_{j=0}^{\tk-1} f \biggl(\frac
{S_k+j}{n},
\frac{\Xkj}{b_n} \biggr) \ind{\frac{\Xko^{(b_n)}}{b_n} < \delta} > 0 \Biggr].
\end{eqnarray*}
As above, the second term is at most $\P[S_{2Rn}/n \leq R ]
\rightarrow0$
as $n\rightarrow\infty$.
The first term is bounded by
\begin{eqnarray*}
&&\P \Biggl[ \bigcup_{k=0}^{2Rn-1} \Biggl(
\biggl\{ \frac{\Xko^{(b_n)}}{b_n} < \delta \biggr\} \cap\bigcup
_{j=1}^{\tk-1} \biggl\{ \frac
{\Xkj}{b_n} \geq M \biggr
\} \Biggr) \Biggr]
\\
&&\quad\leq2Rn \P \biggl[ \frac{X_0^{(b_n)}}{b_n} < \delta,
\sup_{j\geq1} \frac{X_j^{(b_n)}}{b_n} \geq M \biggr]
\\
&&\quad\leq2R n\P \biggl[\frac{X_0^{(b_n)}}{b_n} < \delta, \sup_{1 \leq j
\leq m_0}
\frac{X_j^{(b_n)}}{b_n} \geq M \biggr] + 2R n\P \biggl[\frac
{X_0^{(b_n)}}{b_n} < \delta,
\sup_{j \geq m_0+1} \frac{X_j^{(b_n)}}{b_n} \geq M \biggr]
\\
&&\quad= A_n(\delta) + B_n(\delta),
\end{eqnarray*}
with $m_0$ as in Condition \eqref{condeclrdmrv}.
For $A_n(\delta)$, we have by \eqref{eqecjrvfull},
\[
\lim_{\delta\downarrow0} \limsup_{n\rightarrow\infty} A_n(\delta ) =
\lim_{\delta\downarrow0} 2R\mu^* \bigl( [0,\delta] \times { \bigl([0,M
)^{m_0} \bigr)}^{\mathrm{c}} \bigr) = 0.
\]
For the second, $\lim_{\delta\downarrow0} \limsup_{n\rightarrow
\infty}
B_n(\delta) = 0$ by Condition \eqref{condeclrdmrv}.
This establishes~\eqref{eqconvtogrv}.
\end{pf*}

For completeness, the following lemma notes the continuity of the map
$T_1$ defined in
\eqref{eqpftimech}. See also \cite{whitt:1980,whitt:2002}.

\begin{lem}\label{lemcontvagueconv}
The mapping $T_1 : D^{\uparrow}[0,\infty) \times\mplus([0,\infty)
\times
\E)
\rightarrow\mplus([0,\infty) \times\E)$ given by
\[
T_1{m}(f) = \iint f \bigl(x(u),v \bigr) m(du,dv),
\qquad f\in\mathcal{C}^{+}_K \bigl([0,\infty ) \times
\E\bigr),
\]
is continuous at $(x,m)$ whenever the function $x(\cdot)$ is continuous.
\end{lem}

\begin{pf}
(a) Suppose $x_n \rightarrow x_0$ in $D^{\uparrow}[0,\infty)$ (with
respect to
the Skorohod topology), where $x_0$ is continuous, and $m_n \,\mathop
{\rightarrow}\limits^v\,
m_0$ in $\mplus([0,\infty) \times\E)$.
Let $f\in\mathcal{C}^{+}_K([0,\infty)\times\E)$ with support
contained in
$[0,R]\times B$.
We show that $T_1{m}_n(f) \rightarrow T_1{m}_0(f)$.
For $n\geq0$, write $f_n(u,v) = f(x_n(u),v)$.
The $f_n$ are supported on $x_n^{-1}([0,R])\times B$, and
$x_n^{-1}([0,R]) = [0,x^{\leftarrow}_n(R)]$, where $x^{\leftarrow}_n$
is the
right-continuous inverse of $x_n$.
We now argue that the $f_n$, $n\geq0$, have a common compact support.
Indeed, we have $x_n^{\leftarrow} \rightarrow x_0^{\leftarrow}$
pointwise, so
$x^{\leftarrow}_n(R) \rightarrow x^{\leftarrow}_0(R)$. Thus, for
large $n$, $[0,x^{\leftarrow}_n(R)] \times B \subset[0,x^{\leftarrow
}_0(R)+1] \times B$; without loss
of generality, $m_0(\bdry([0,x^{\leftarrow}_0(R)+1] \times B))= 0$.
Furthermore, $f_n \rightarrow f_0$ uniformly: suppose $(u_n,v_n)
\rightarrow
(u_0,v_0) \in[0,\infty) \times\E$. Then $x_n(u_n) \rightarrow x_0(u_0)$
since $x_0$ is continuous, and so $f(x_n(u_n),v_n) \rightarrow
f(x_0(u_0),v_0)$ by the continuity of $f$.
Consequently, $\tilde{m}_n(f) \rightarrow\tilde{m}_0(f)$ by \cite
[Lemma~8.2(b)]{resnick2011asymptotics}.
\end{pf}




%
%
%


\section*{Acknowledgements}
S.I. Resnick and D. Zeber were partially supported by ARO Contract
W911NF-10-1-0289 and NSA Grant H98230-11-1-0193 at Cornell University.

%
%

\end{document}